\documentclass[5p]{elsarticle}

\usepackage[nolist]{acronym}
\usepackage{hyperref}
\usepackage{amsfonts}
\usepackage{amssymb}  
\usepackage{graphicx,epstopdf}
\usepackage[caption=false]{subfig}
\usepackage{algorithmic}
\usepackage{acronym}
\usepackage{amsmath}
\usepackage{natbib}  % biblio  
\usepackage{leftidx} 

\journal{Journal of Control Engineering Practice}

\begin{acronym}
	\acro{UAV}{Unmanned Aerial Vehicle}
	\acro{CFD}{Computational Fluid Dynamics}
	\acro{AWE}{Airborne Wind Energy}
	\acro{OED}{Optimum Experimental Design}
	\acro{CAD}{Computer-Aided Design}
	\acro{IMU}{Inertial Measurement Unit}
	\acro{FCC}{Flight Control Computer}
	\acro{ODE}{Ordinary Differential Equation}
	\acro{ODEs}{Ordinary Differential Equations}
	\acro{NLP}{Non-Linear Program}
	\acro{CRLB}{Cramer-Rao Lower Bound}
	\acro{LTI}{Linear Time-Invariant}
	\acro{SNR}{Signal-to-Noise Ratio}
	\acro{EOM}{Equation of Motion}
	\acro{DAE}{Differential Algebraic Equation}
	\acro{DAEs}{Differential Algebraic Equations}
	\acro{AWES}{Airborne Wind Energy System}
	\acro{AWE}{Airborne Wind Energy}
	\acro{NED}{North-East-Down}
	\acro{DCM}{Direct Cosine Matrix}
	\acro{MBPE}{Model-Based Parameter Estimation}
	\acro{OCP}{Optimal Control Problem}
	\acro{IVP}{Initial Value Problem}
	\acro{TIC}{Theil Inequality Coefficient}
	\acro{PE}{Parameter Estimation}
\end{acronym}

\newcommand{\VT}{V_\mathrm{T}}
\newcommand{\VTe}{V_{\mathrm{T}_\mathrm{e}}}
\begin{document}
\begin{frontmatter}
\title{System Identification of a Rigid Wing Airborne Wind Energy System\tnoteref{acknowledgment}}
\tnotetext[acknowledgment]{This research was supported by the EU via FP7-ITN-TEMPO (607 957) and H2020-ITN-AWESCO (642 682), by the Federal Ministry for Economic Affairs and Energy (BMWi) via eco4wind and DyConPV, and by DFG via Research Unit FOR 2401.}

\author[AP,UniFr]{G.~Licitra\corref{CorrespondingAuthor}}\ead{giovanni.licitra@imtek.uni-freiburg.de}
\author[UniKal,UniFr]{A.~B\"{u}rger}\ead{adrian.buerger@hs-karlsruhe.de}
\author[AP]{P.~Williams}\ead{p.williams@ampyxpower.com}
\author[AP]{R.~Ruiterkamp}\ead{r.ruiterkamp@ampyxpower.com}
\author[UniFr]{M.~Diehl}\ead{moritz.diehl@imtek.uni-freiburg.de}

\cortext[CorrespondingAuthor]{Corresponding author}

\address[AP]{Ampyx Power B.V, The Hague, Netherlands}
\address[UniFr]{Department of Microsystems Engineering, University of Freiburg, Germany}
\address[UniKal]{Faculty of Management Science and Engineering, Karlsruhe University of Applied Sciences, Germany}

\begin{abstract}
	\ac{AWE} refers to a novel technology capable of harvesting energy from wind by flying crosswind patterns with tethered autonomous aircraft. 
	Successful design of flight controllers for \ac{AWE} systems rely on the availability of accurate mathematical models. Due to the non-conventional structure of the airborne component, the system identification procedure must be ultimately addressed via an intensive flight test campaign to gain additional insight about the aerodynamic properties.  
	In this paper, aerodynamic coefficients are estimated from experimental data obtained within flight tests using an multiple experiments \ac{MBPE} algorithm. 
\end{abstract}
\begin{keyword}
	Airborne Wind Energy, autonomous aircraft, optimization, system identification and parameter estimation.
\end{keyword}
\end{frontmatter}

\section{Introduction}
In the landscape of innovative renewable energy systems, Airborne Wind Energy (AWE) is a novel, emerging technology. \ac{AWE} promises to harvest energy from wind with both lower installation costs and higher capacity factors compared to conventional wind turbines, up to a level that could render \ac{AWE} even more economically viable than fossil fuels.

Despite the fact that the idea of using tethered aircraft for wind power generation appeared for the first time in the late 1970s \cite{loyd1980crosswind}, it is only in the last decade that academia and industry made substantial progress in turning the idea into a practical implementation.
The postponement of \ac{AWE} technology is mainly due to the significant complexities in terms of control, modeling, identification, materials, mechanics and power electronics. 
Furthermore, these systems need to fulfill high level of reliability while simultaneously operating close to optimality. Such requirements have brought many developers to the use of rigid wing autonomous aircraft as airborne component \cite{AP,makani,twingtec,kitemill,e-kite}.  

In the aerospace field, it is the current practice to retrieve the aircraft aerodynamic properties by a combination of wind tunnel testing, \ac{CFD} \cite{anderson2017fundamentals} analysis, and empirical methods such as DATCOM \cite{hoak1975usaf}. For standard aircraft configurations, such methods for obtaining aerodynamic characteristics are generally in agreement with those obtained via flight test. However, empirical methods, which can provide the quickest results, tend to be less accurate and more difficult to apply to unconventional designs. \ac{CFD} is much more accurate, but requires a fine mesh to capture the flow dynamics accurately, and as a consequence it involves significant computational resources to obtain a complete aerodynamic database.  As far as it regards wind tunnel experiments, they generally provide the most accurate results with a suitably sized model that matches the Reynold's numbers of the real system.  However, for unconventional systems such approach can also be expensive. In any case, an intensive flight test campaign must be set in order to gain additional insight about the aerodynamic properties \cite{licitra2017pe} and to validate parameters on the full scale system.

\begin{figure}[t]
	\centering
	\includegraphics[scale = 0.42]{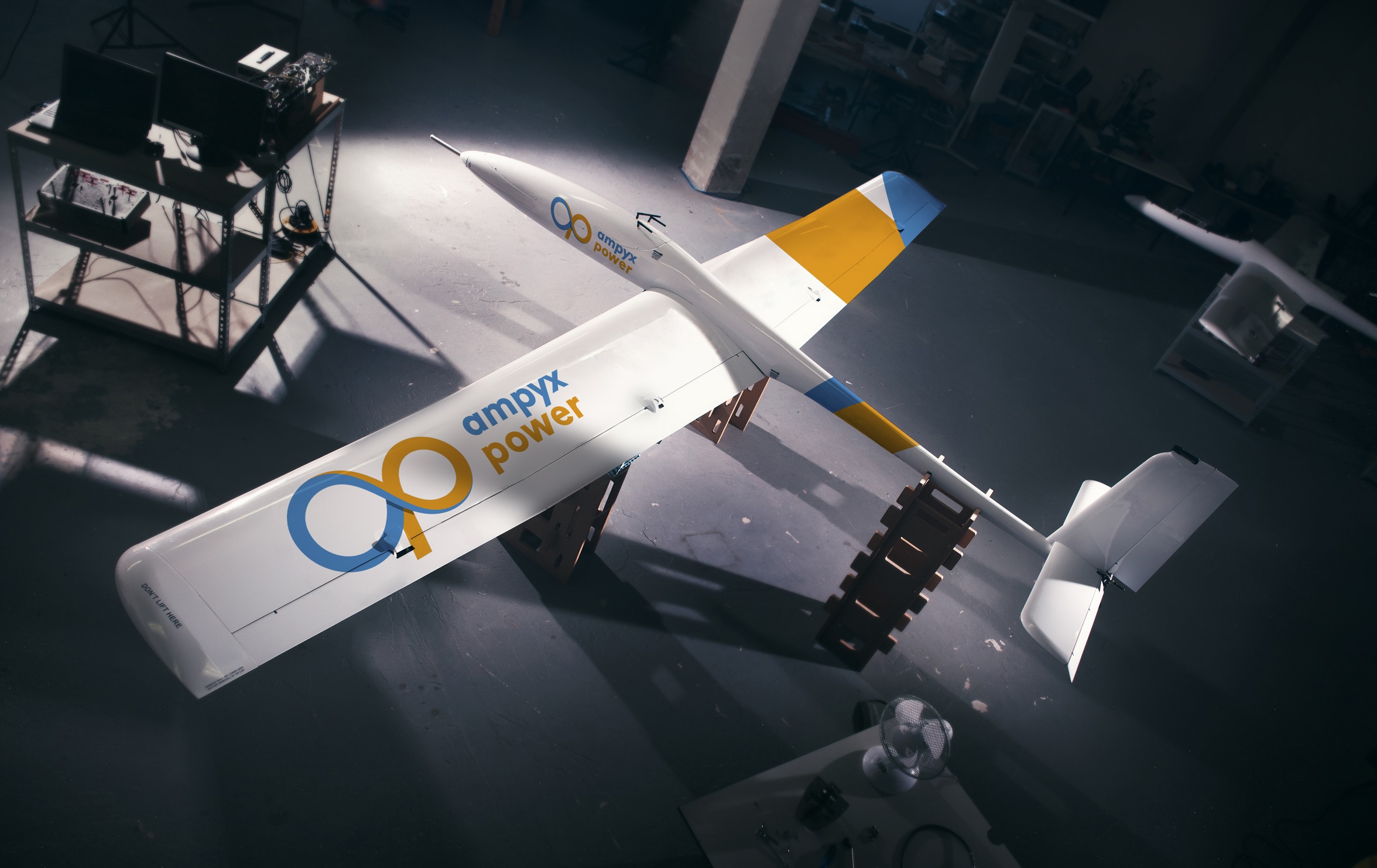}
	\caption{The $2^{\mathrm{nd}}$ prototype high lift, rigid wing autonomous aircraft designed by Ampyx Power B.V.}
	\label{fig:AP2}
\end{figure}

The presented work is entirely based on the $2^{\mathrm{nd}}$ prototype high lift, rigid wing autonomous aircraft designed by Ampyx Power B.V. \cite{AP} and shown in figure~\ref{fig:AP2}.

\begin{figure}[tphb]
	\begin{center}
		\includegraphics[scale=0.38]{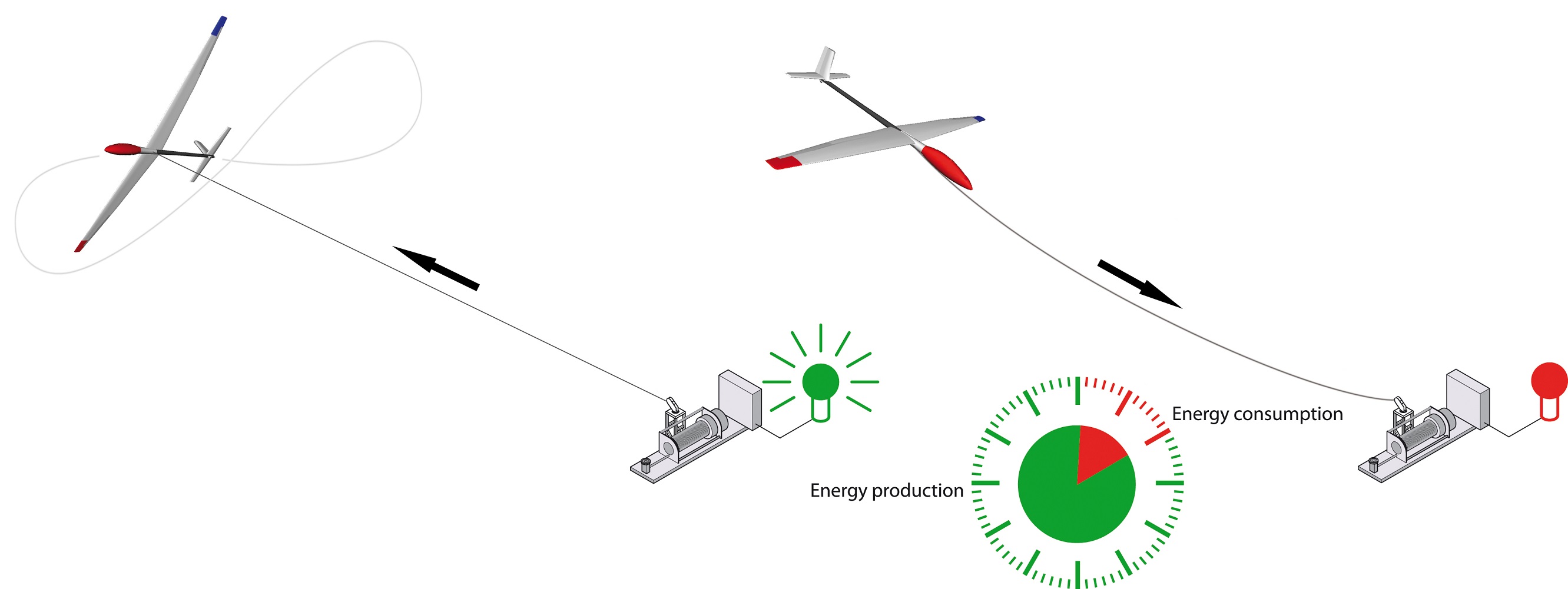}  
		\caption{Working principles of a lift mode \ac{AWES} with a \textit{production} and \textit{consumption phase}. A lift mode \ac{AWES} produces power by performing periodical variation of both length and tether tension. Power generation occurs during the so called \textit{reel-out phase}, where the tether tension is used to rotate a drum, driving an electric generator located on the ground. A \textit{reel-in phase} is required due to finite tether length. By changing the flight pattern in such a way that less lifting force is produced, the tether can be wound with a significant lower energy investment than what was gained in the power production phase.}
		\label{fig:PumpingModeConcept}
	\end{center}
\end{figure}

Ampyx Power B.V. adopts the so called \textit{lift mode} strategy \cite{loyd1980crosswind,diehl2013airborne,Cap26AWEbook,PowerPatternvideo} where the airplane delivers a high tension on the tether which is anchored to a ground-based generator. An artist's rendering of the two main phases of a lift mode \ac{AWES} is shown in figure~\ref{fig:PumpingModeConcept}. 

\begin{figure}[tphb]
	\begin{center}
		\includegraphics[scale=0.11]{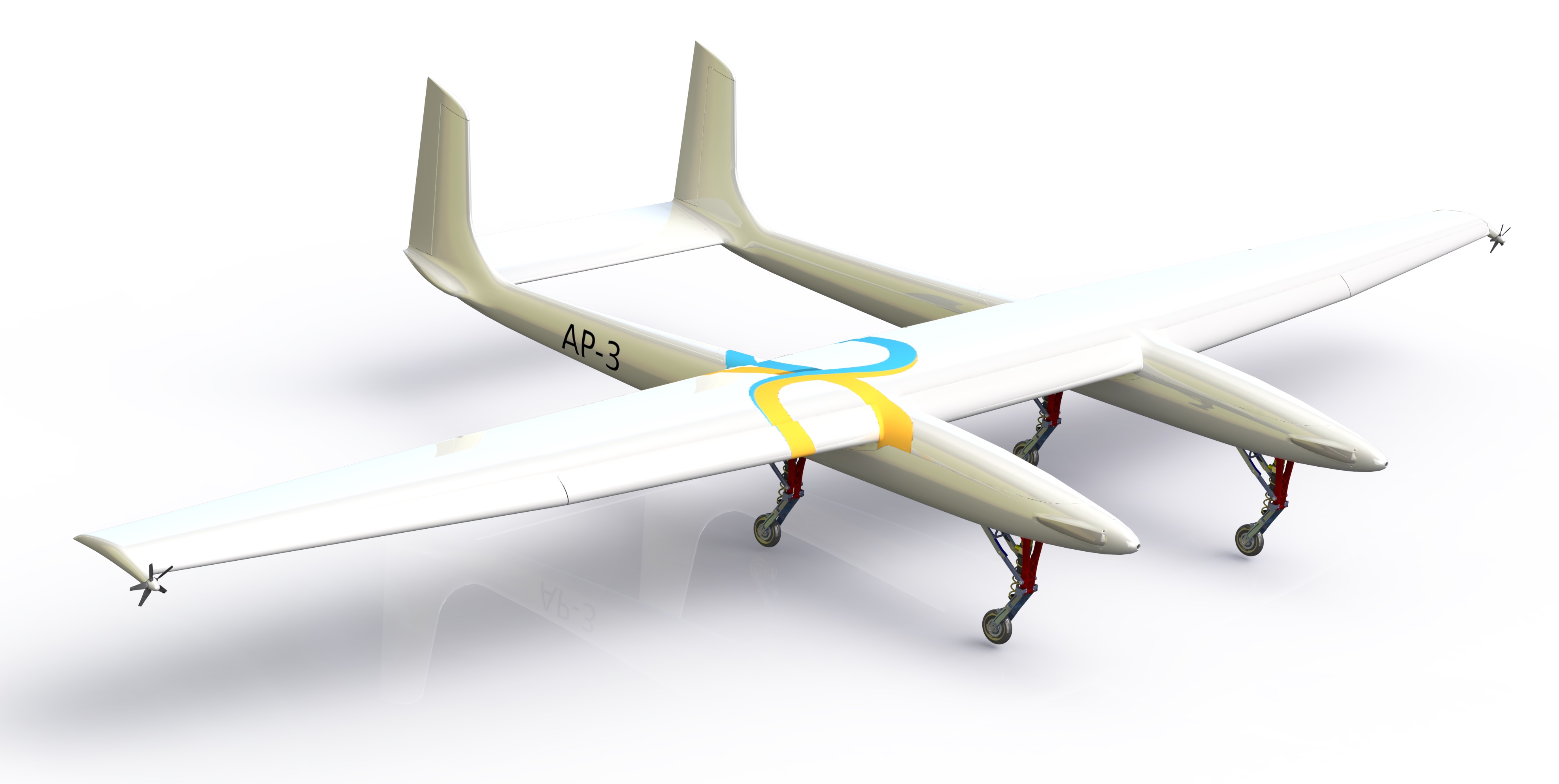}%[width = 240pt, height = 120pt]  
		\caption{The $3^{rd}$ prototype high lift, rigid wing autonomous aircraft designed by Ampyx Power B.V.}
		\label{fig:AP3}
	\end{center}
\end{figure}

A successful flight test campaign which aims to identify the aerodynamic parameters of the aircraft depends on many factors, such as selection of instrumentation, signal conditioning, flight test operations procedure, parameter estimation algorithm and signal input design. In \cite{licitra2017pe}, aerodynamic properties were estimated via flight tests with conventional maneuvers for the pitch rate dynamics, only. In \cite{licitra2017oed}, optimal maneuvers were computed for the case study by solving a time domain model-based \ac{OED} problem which aims to obtain more accurate parameter estimates while enforcing safety constraints. The optimized inputs were compared with respect to conventional maneuvers widely used in the aerospace field and successfully tested within real experiments \cite{licitra2017optinput}. 
In this paper, estimation of the aerodynamic characteristics are carried out via an efficient multiple experiment Model-Based Parameter Estimation (\ac{MBPE}) algorithm based on \textit{direct methods} using both conventional and optimized experiments. The data fitting is applied throughout the aircraft longitudinal dynamics using a non-linear model structure.   
The presented work will be used as guideline for the system identification of the next prototype designed by Ampyx Power B.V. \cite{AP} and shown in figure~\ref{fig:AP3}. 

The paper is organized as follows. In section \ref{sec:IntroModel} the mathematical model of a rigid wing airborne component of a generic \ac{AWES} is introduced. Subsequently, a suitable model structure is selected for the estimation of aerodynamic properties augmented with model assumptions as well as neglected dynamics. Section \ref{sec:DesignEvaluationInputSignals} focuses on the design and evaluation of input signals. A preliminary estimation performance analysis is carried out using the Fisher information matrix and an overview of both flight test procedures and decoupling of dynamics are provided. In section \ref{sec:ExperimentalData} the experimental data obtained from conventional and optimized flight test campaigns are shown. Section \ref{sec:MBPEformulation} formulates the multiple experiment \ac{MBPE} algorithm whereas in section \ref{sec:dataFitting} the data fitting is performed on the obtained experimental data. Finally, in section \ref{sec:ModelValidation} both estimates and model validation are assessed and conclusions are provided in section \ref{sec:Conclusions}. 

\section{Modeling of a rigid wing \ac{AWE} system}\label{sec:IntroModel}
In this section, a mathematical formulation of an \ac{AWES} is introduced. Subsequently, a non-linear model structure is selected for the purposes of system identification underlying model assumptions and neglected dynamics.

\subsection{Modeling of AWES in natural coordinates}\label{sec:MathematicalModel}
A rigid wing \ac{AWES} can be efficiently modeled as a set of \ac{DAEs} described by non-minimal coordinates by means of Lagrangian mechanics. The equations of motion for a tethered airborne component are given by \cite{gros2013modeling}:
\begin{subequations}\label{eq:EOMfull}
	\begin{align}
	\mathbf{\dot{p}^{n}} & = \mathbf{R_{nb}} \cdot \mathbf{v^{b}} \label{eq:pNED}\\ 
	m \cdot \mathbf{\dot{v}^{b}} & = \mathbf{f_{c}^{b}} + \mathbf{f_{p}^{b}} + \mathbf{f_{a}^{b}} + \mathbf{f_{g}^{b}} - m \, ( \mathbf{\omega^{b}} \times \mathbf{v^{b}}) \label{eq:transEOMfull}  \\
	\mathbf{\dot{R}_{nb}} & = \mathbf{R_{nb}} \cdot \mathbf{\Omega} \label{eq:DCM}\\
	\mathbf{J} \cdot \mathbf{\dot{\omega}^{b}} & = \mathbf{m_{c}^{b}} + \mathbf{m_{p}^{b}} + \mathbf{m_{a}^{b}} - (\mathbf{\omega^{b}} \times \mathbf{J} \cdot \mathbf{\omega^{b}} ) \label{eq:rotEOMfull}
	\end{align}
\end{subequations}
where $\mathbf{v^{b}}=[u,v,w]^{\top} $ and $\mathbf{\omega^{b}} = [p,q,r]^{\top}$ are the translational and rotational speed vector defined in body-fixed frame $\mathbf{b}$, $m$ the mass and $\mathbf{J}$ the inertia dyadic of the aircraft. In \eqref{eq:pNED}, the rate of change in position $\mathbf{\dot{p}^{n}}$ is defined in \ac{NED} frame and it is obtained by means of the \ac{DCM} from body to \ac{NED} frame $\mathbf{R_{nb}} \in \mathbb{R}^{3 \times 3}$ whereas \eqref{eq:DCM} is the time evolution of the \ac{DCM} with $\mathbf{\Omega}  \in \mathbb{R}^{3 \times 3}$ the skew symmetric matrix of $\mathbf{\omega^{b}}$.
The aircraft is subject to forces $\mathbf{f_{\{c,p,g\}}^{b}}$ and moments $\mathbf{m_{\{c,p,g\}}^{b}}$ coming from the cable, propellers, gravity and the interaction between aircraft with the air mass is denoted by $\mathbf{f_{a}^{b}} = \left[X,Y,Z\right]^{T} $ and $\mathbf{m_{a}^{b}} = \left[L,M,N\right]^{T}$. The mathematical formulation in \eqref{eq:EOMfull} is extensively used for pattern generation using an optimal control approach \cite{licitra2016optimal,licitra2017viability}.  

In order to identify the aerodynamic forces $\mathbf{f_{a}^{b}}$ and moments $\mathbf{m_{a}^{b}}$, one has to either discard or have good models of the other contributions. For this application, it is both hardly possible have accurate cable models and vibration effects coming from the cable itself arise during tethered crosswind flights. Furthermore, the propulsion system introduce additional noise for each angular rate and acceleration channel provided by the rotation of the blades. 

Hence, a flight test campaign which aims to the identification of aerodynamic properties needs to be performed without tether such that the cable does not interfere with the overall aircraft dynamics \cite{licitra2017pe} and additionally, propellers must be switched off whenever an excitation signal occurs so as to decouple the uncertainty in thrust effects on the aerodynamic parameter estimation, simplifying \eqref{eq:EOMfull} to
\begin{subequations}\label{eq:EOMsimplified}
	\begin{align}
	m \cdot \mathbf{\dot{v}^{b}} & = \mathbf{f_{a}^{b}} + \mathbf{f_{g}^{b}} - m \, ( \mathbf{\omega^{b}} \times \mathbf{v^{b}}) \label{eq:transEOMsimplified}  \\
	\mathbf{\dot{R}_{nb}} & = \mathbf{R_{nb}} \cdot \mathbf{\Omega}\\
	\mathbf{J} \cdot \mathbf{\dot{\omega}^{b}} & =  \mathbf{m_{a}^{b}} - (\mathbf{\omega^{b}} \times \mathbf{J} \cdot \mathbf{\omega^{b}} ) \label{eq:rotEOM_simplified}
	\end{align}
\end{subequations}
Note that equation~\eqref{eq:pNED} is discarded since it does not provide any meaningful information for system identification purposes. 

Finally, as far as it regards the cable dynamics, a comprehensive study can be found in \cite{williams2007modeling,williams2008modeling} whereas the propeller forces and moments are normally obtained via extensive test bench. 

\subsection{Model Selection}\label{sec:ModelSelection}
The case study considered within this work is a high lift, rigid wing autonomous aircraft used as airborne component of a lift mode \ac{AWES} designed by Ampyx Power B.V. \cite{AP}. Details on the system can be found in \cite{AWEbook,Cap26AWEbook,licitra2017optinput} whereas Table~\ref{tab:AP2_parameters} collects the main physical properties.

\begin{figure}[tbhp]
	\centering
	\includegraphics[scale = 0.25]{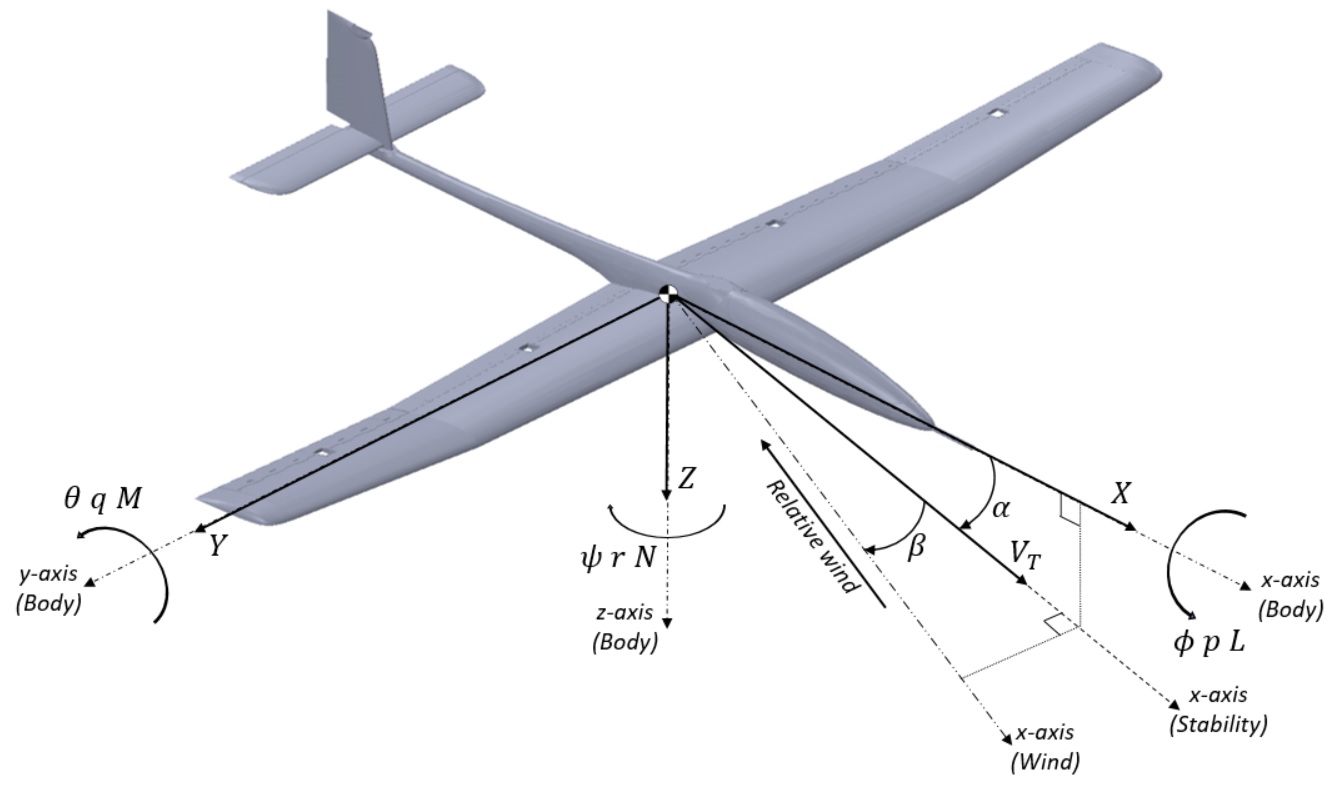}% [width = 240pt, height = 100pt] 
	\caption{Definition of axes, Euler angles, aerodynamic states, forces and moments on a rigid wing aircraft.}
	\label{fig:AircraftConvention}
\end{figure}

For system identification purposes, it is more convenient to have the velocity equation \eqref{eq:transEOMsimplified} in terms of \textit{wind-axes} variables: airspeed $\VT$ and aerodynamic angles $\beta$ and $\alpha$ which are the angle of side-slip and attack, respectively. Furthermore, the aircraft attitude can be described via the Euler angles kinematics where $\phi$,$\theta$,$\psi$ denote the roll, pitch and yaw angle.   
The proposed model structure is therefore given by \cite{stevens2015aircraft}:
\begin{subequations}\label{eq:MathematicalModel}
	\begin{align}
	\dot{V}_{\mathrm{T}}  & = \frac{Y s\beta + X c\alpha c\beta + Z c\beta  s\alpha}{m} + G_{\VT}, \label{eq:Vt}\\
	\dot{\beta}  & = \frac{Y \mathrm{c}\beta  - X \mathrm{c}\alpha \mathrm{s}\beta - Z \mathrm{s}\alpha \mathrm{s}\beta}{m \VT} + \frac{G_{\beta}}{\VT} - r \mathrm{c}\alpha + p \mathrm{s}\alpha, \label{eq:beta}\\
	\dot{\alpha} & =  \frac{Z \mathrm{c}\alpha - X \mathrm{s}\alpha}{m \VT \mathrm{c}\beta} + \frac{G_{\alpha}}{\VT \mathrm{c}\beta} + \frac{q \mathrm{c}\beta - (p \mathrm{c}\alpha + r \mathrm{s}\alpha) \mathrm{s}\beta}{\mathrm{c}\beta}, \label{eq:alpha}\\
	\dot{\phi}   & = p + r \, \mathrm{c}\phi \, \mathrm{t}\theta + q \, \mathrm{s}\phi \, \mathrm{t}\theta, \label{eq:phi}\\
	\dot{\theta} & = q \, \mathrm{c}\phi - r \, \mathrm{s}\phi, \label{eq:theta}\\
	\dot{\psi}   & = \frac{q \, \mathrm{s}\phi +  r \, \mathrm{c}\phi}{\mathrm{c}\theta}, \label{eq:psi}\\		
	\dot{p}      & = \frac{J_{xz}}{J_{x}} \dot{r} - qr\frac{\left(J_{z} - J_{y} \right)}{J_{x}} + qp \frac{J_{xz}}{J_{x}} + \frac{L}{J_{x}}, \label{eq:p}\\
	\dot{q}      & =  -pr \frac{J_{x} - J_{z}}{J_{y}} - (p^{2} - r^{2})\frac{J_{xz}}{J_{y}} + \frac{M}{J_{y}}, \label{eq:q}\\
	\dot{r}      & = \frac{J_{xz}}{J_{z}} \dot{p} -pq \frac{J_{y}-J_{x}}{J_{z}} - qr\frac{J_{xz}}{J_{z}} + \frac{N}{J_{z}}, \label{eq:r}	
	\end{align}
\end{subequations}
where $\mathrm{s(.)},\mathrm{c(.)},\mathrm{t(.)}$ are a shortening of the trigonometric functions $\sin(.),\cos(.),\tan(.)$ and $G_{\VT}$,$G_{\beta}$,$G_{\alpha}$ are the gravity components expressed in wind frame and equal to
\begin{subequations}\label{eq:GravityComponents}
	\begin{align}	    
	G_{\VT}    & =  g_{\mathrm{D}} \left(\mathrm{s}\beta   \, \mathrm{s}\phi \, \mathrm{s}\theta -\mathrm{c}\alpha  \, \mathrm{c}\beta \, \mathrm{s}\theta + \mathrm{s}\alpha \, \mathrm{c}\beta \, \mathrm{c}\phi  \, \mathrm{c}\theta \right), \\
	G_{\beta}  & =  g_{\mathrm{D}} \left(\mathrm{c}\alpha  \, \mathrm{s}\beta \, \mathrm{s}\theta + \mathrm{c}\beta  \, \mathrm{s}\phi  \, \mathrm{c}\theta - \mathrm{s}\alpha \, \mathrm{s}\beta \, \mathrm{c}\phi \, \mathrm{c}\theta \right),\\
	G_{\alpha} & =  g_{\mathrm{D}} \left(\mathrm{s}\alpha \, \mathrm{s}\theta + \mathrm{c}\alpha \, \mathrm{c}\phi \, \mathrm{c}\theta \right),
	\end{align}
\end{subequations}
with $g_{\mathrm{D}} \approx 9.81\,\mathrm{m/s^{2}}$ the gravitational acceleration. The nomenclature introduced above is summarized in figure~\ref{fig:AircraftConvention}. The mathematical model \eqref{eq:MathematicalModel} implicitly presumes that the vehicle is a rigid body with a plane of symmetry such that the moments of inertia $J_{xy},J_{xz}$ are zero, whereas the Earth is assumed flat and non-rotating with a constant gravity field \cite{mulder2000flight}.
Note that, the model equations in \eqref{eq:MathematicalModel} is also widely used for linearization purposes, dynamics analysis as well as control system design \cite{stevens2015aircraft}.  
Within this work, the aerodynamic forces $\left(X,Y,Z\right)$ and moments $\left(L,M,N\right)$ are normalized with respect to the dynamic pressure $\bar{q} = \frac{1}{2} \rho \VT^{2}$ with $\rho\approx 1.225 \,\mathrm{kg/m^{3}}$ the free-stream mass density, and a characteristic area for the aircraft body
\begin{subequations}\label{eq:Fa_Ma}
	\begin{align}
	X & = \bar{q}S        \,C_{X} & Y & = \bar{q}S        \,C_{Y} & Z & = \bar{q}S        \,C_{Z} \\
	L & = \bar{q}Sb       \,C_{l} & M & = \bar{q}S\bar{c} \,C_{m} & N & = \bar{q}Sb       \,C_{n}. 				
	\end{align}
\end{subequations}
In \eqref{eq:Fa_Ma} $S,\,b,\,\bar{c}$ are reference wing area, wing span and mean aerodynamic chord, respectively, while $C_{X},\,C_{Y},\,C_{Z}$ denote the forces and $C_{l},\,C_{m},\,C_{n}$ the moment coefficients.
For conventional aircraft, the aerodynamic coefficients are usually broken down into a sum of terms as follows
\begin{subequations}\label{eq:AdimensionaAerodynamicCoefficients}
	\begin{align}
	C_{X} & = C_{X_{\alpha}}\alpha + C_{X_{q}}\hat{q} + C_{X_{\delta_{\mathrm{e}}}} \delta_{\mathrm{e}} + C_{X_{0}},\\
	C_{Y} & = C_{Y_{\beta}} \beta  + C_{Y_{p}}\hat{p} + C_{Y_{r}} \hat{r} + C_{Y_{\delta_{\mathrm{a}}}} \delta_{\mathrm{a}} + C_{Y_{\delta_{\mathrm{r}}}} \delta_{\mathrm{r}}, \\  
	C_{Z} & = C_{Z_{\alpha}}\alpha + C_{Z_{q}}\hat{q} + C_{Z_{\delta_{\mathrm{e}}}} \delta_{\mathrm{e}} + C_{Z_{0}},\\     	
	C_{l} & = C_{l_{\beta}} \beta  + C_{l_{p}}\hat{p} + C_{l_{r}} \hat{r} + C_{l_{\delta_{\mathrm{a}}}} \delta_{\mathrm{a}} + C_{l_{\delta_{\mathrm{r}}}} \delta_{\mathrm{r}}, \label{eq:Cl}\\
	C_{m} & = C_{m_{\alpha}}\alpha + C_{m_{q}}\hat{q} + C_{m_{\delta_{\mathrm{e}}}} \delta_{\mathrm{e}} + C_{m_{0}}, \\
	C_{n} & = C_{n_{\beta}} \beta  + C_{n_{p}}\hat{p} + C_{n_{r}} \hat{r} + C_{n_{\delta_{\mathrm{a}}}} \delta_{\mathrm{a}} + C_{n_{\delta_{\mathrm{r}}}} \delta_{\mathrm{r}},
	\end{align}
\end{subequations}
that depend on the normalized body rates $\hat{p}= \frac{b \, p}{2\VT},\hat{q} = \frac{\bar{c} \, q}{2\VT},\hat{r} = \frac{b \, r}{2\VT}$, angle of attack $\alpha$ and side slip $\beta$, as well as the control surface deflections which in this case are aileron $\delta_{\mathrm{a}}$, elevator $\delta_{\mathrm{e}}$ and rudder $\delta_{\mathrm{r}}$. The coefficients $C_{i_{j}}$ with $i = \{X,Y,Z,l,m,n\}$ and $j = \{\alpha,\beta,p,q,r,\delta_{\mathrm{a}},\delta_{\mathrm{e}},\delta_{\mathrm{r}},0\} $ are the \textit{dimensionless aerodynamic derivatives} that need to be identified.
\begin{table}[tbhp]
	\caption{Aircraft physical properties}
	\label{tab:AP2_parameters}
	\centering
	\begin{tabular}{lccc}\hline
		\bf Name                &\bf Symbol &\bf Value & \bf Unit                  \\\hline
		mass                    & $ m $     & $36.8$   & $\mathrm{kg}$             \\
		moment of inertia       & $J_{x}$   & $25$     & $\mathrm{kg \cdot m^{2}}$ \\
		moment of inertia       & $J_{y}$   & $32$     & $\mathrm{kg \cdot m^{2}}$ \\
		moment of inertia       & $J_{z}$   & $56$     & $\mathrm{kg \cdot m^{2}}$ \\
		cross moment of inertia & $J_{xz}$  & $0.47$   & $\mathrm{kg \cdot m^{2}}$ \\
		reference wing area     & $S$       & $3$      & $\mathrm{m^{2}}$          \\
		reference wing span     & $b$       & $5.5$    & $\mathrm{m}$              \\
		reference chord         & $\bar{c}$ & $0.55$   & $\mathrm{m}$              \\\hline		
	\end{tabular}
\end{table}

\subsection{Model assumption and neglected dynamics}\label{sec:ModelAssumption}
In general, the aerodynamic forces and moments are all dependent on the time history of the aircraft state in time, which means that if the pitch moment $M$ depends on the pitch rate $q$ only, then:
\begin{equation}
M(t) = f(q(t)), t \in (-\infty,\tau].
\end{equation}
Theoretically speaking, the function in time $q(t)$ can be replaced by the following Taylor series:
\begin{equation}\label{eq:DerivativesInTime}
q(t) = q(\tau) + \sum_{i=1}^{\infty} \frac{1}{i!} \frac{\partial^i q }{\partial \tau^i} (t - \tau)^{i}  
\end{equation}
i.e. that the whole information regarding the parameter history $q$ is captured, if we were able to compute all the possible derivatives. However, for subsonic flight the influence of the derivatives is bounded and can be neglected \cite{mulder2000flight}.

In flight dynamics there are different methods of aerodynamic derivatives modeling. In many practical cases, the aerodynamic properties are approximated by linear terms in their Taylor series expansion as in \eqref{eq:AdimensionaAerodynamicCoefficients}. On the one hand, such approximations yield sufficient accuracy for attached flows \cite{etkin2012dynamics}. On the other hand, this representation cannot be used in the region of $\alpha$ where separated and vortex flow occur \cite{goman1994state}.   

\begin{figure}[tbhp]
	\centering
	\includegraphics[width = 260pt, height = 150pt]{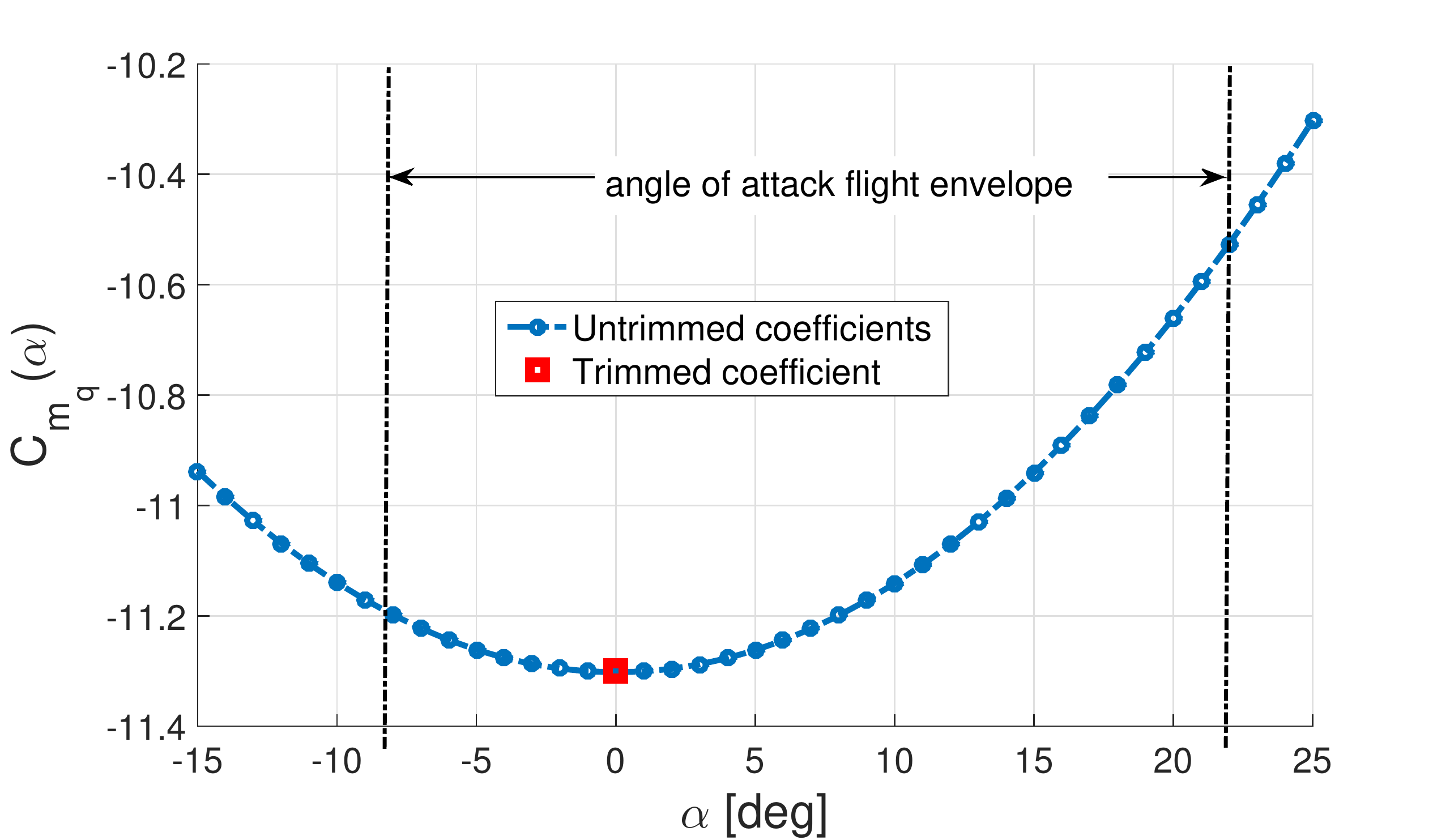}	%[scale=0.28]
	\caption{A priori pitch damping derivative $C_{m_{q}}(\alpha)$ with corresponding trimmed coefficient for $\VTe = 20$ \rm{[m/s]}.}
	\label{fig:Cmq}
\end{figure}

In this work, since the aircraft dynamics and its aerodynamic characteristics are described via the equations  \eqref{eq:MathematicalModel}, \eqref{eq:GravityComponents}, \eqref{eq:Fa_Ma} and \eqref{eq:AdimensionaAerodynamicCoefficients}, one has to implicitly account for the model mismatches summarized below
\begin{itemize}
	\item By neglecting the influence of the derivatives over time shown in \eqref{eq:DerivativesInTime}, one neglects the influence of parameter variation through time. Such influence arises from non-stationary wing-fuselage and tail interference, increasing during aggressive maneuvers \cite{mulder2000flight}, in our case mainly during the \textit{power-generation} phase. However, some dynamics can be captured by introducing a first-order differential equation involving the angle of attack rate $\dot{\alpha}$, c.f. \cite{goman1994state}.  
	\item The case study utilizes a high-strength wing with relatively high stiffness. As a consequence, flexible modes might occur on the airframe, though the mathematical model in \eqref{eq:EOMfull} does not take in account such effects. Flexible modes need to be considered during the control systems design phase so as to avoid possible structural-coupling issues. 
	\item The aerodynamic model assumed in \eqref{eq:AdimensionaAerodynamicCoefficients} is implicitly a function of $\alpha$. Nevertheless, system identification performed via flight tests are typically valid only for small neighborhood of $\alpha$ with respect to its trim value $\alpha_{\mathrm{e}}$ given at a specific trim airspeed $\VTe$. Because aircraft deployed for \ac{AWES} are intended to fly over a wide range of flight conditions, flight test maneuvers and parameter identification needs to be performed at multiple trim conditions. Figure~\ref{fig:Cmq} shows the a priori pitch damping coefficient $C_{m_{q}}$ related to the case study as a function of $\alpha$ and \textit{trimmed} value for $\VTe = 20$ \rm{[m/s]}. 
	\item Estimates of aerodynamic derivatives are computed assuming that the aircraft inertias are known a priori. However, fully accurate inertial estimates are difficult to obtain. Inertia estimates can be computed from Computer Aided Design (CAD) models or swing tests with varying degrees of accuracy \cite{de1987accurate,lyons2002obtaining}. Errors in $J_{\{x,y,z,xz\}}$ will lead to errors in the absolute estimates of the aerodynamic coefficients. Nevertheless, this will not undermine the predictive capability of the derived model, as long as the estimated derivatives are kept consistent with the assumed value of $J_{\{x,y,z,xz\}}$ used to estimate them \cite{licitra2017pe}. 
\end{itemize}
In order to overcome the issues mentioned above, it is current practice to design a complex hierarchical control system with high margin or robustness and fly patterns with specific boundary conditions (For further details, refer to \cite{Cap26AWEbook}).  

\section{Design and evaluation of input signals}
\label{sec:DesignEvaluationInputSignals}
In this section, an insight about the flight operation procedure and the rationale behind is provided. Subsequently conventional and optimized maneuvers are designed for parameter estimation purposes and assessed via the \ac{CRLB}. Finally, the experimental data obtained within two flight test campaigns are shown.
\subsection{Flight test procedure and decoupling of dynamics}
\label{sec:steadyConditionAndDecoupling}
Typically, experiments are repeated on each axis to both obtain a rich data set and reduce the effect of sensor biases as well as colored noise (atmospheric turbulence) on the estimation results \cite{licitra2017pe}.
To prevent biases due to correlation between the measurement noise and the inputs, it is best to perform open-loop experiments \cite{soderstrom2001identification}. 
For both physical and practical reasons, system identification flight test are performed at \textit{steady wing-level flight} condition \cite{dorobantu2013system}. More precisely, an aircraft is meant in steady wing-level flight when its body angular rates ($p,q,r$) and roll angle $\phi$ are equal to zero and it flies with constant airspeed $\VTe$ \cite{stevens2015aircraft}. Fulfillment of this steady condition allows decoupling of the aircraft motion in \textit{longitudinal} and \textit{lateral} dynamics, hence one can focus only on a subset of the entire aircraft dynamics which is mainly excited from a given maneuver. For instance, if a signal excitation is performed along the longitudinal axis via elevator deflection, the remaining control surfaces (aileron and rudder deflection) are fixed at their trim values throughout the experiment in order to stabilize the lateral dynamics. As a consequence, the parameter estimation will be performed only on the excited dynamics, which is for this work the longitudinal motion. 

For conventional aircraft parameter estimation experiments, typically a linear perturbation model structure is assumed \cite{klein1990optimal}. Therefore, the flight test inputs are perturbations with respect to the steady condition. In this work, data fitting will be performed using the non-linear formulation \eqref{eq:MathematicalModel} relative to the longitudinal dynamics, however linear representations are used for input design as well as assessment of the expected estimation performance.  

The longitudinal dynamics are described via LTI state-space form by the states $\mathbf{x_{lon}} = \left[V_{\mathrm{T}}~\alpha~\theta~q\right]^\top$, which correspond to \eqref{eq:Vt}, \eqref{eq:alpha}, \eqref{eq:theta} and \eqref{eq:q}. The forces $X$, $Z$ and the moment $M$ are assumed to be linear functions of $\VT,\alpha,q$ and the elevator deflection $\delta_{\mathrm{e}}$, resulting in the following matrices 
\begin{subequations}\label{eq:longDyn}
	\begin{align}
	\mathbf{A_{lon}} & = 
	\begin{bmatrix}
	X_{V} & X_{\alpha}                   & -g_{\mathrm{D}}\cos\theta_{e} & X_{q} \\
	Z_{V} & \frac{Z_{\alpha}}{V_{T_{e}}} & -g_{\mathrm{D}}\sin\theta_{e} & Z_{q} \\
	0   &     0                        &           0          &  1    \\
	M_{V} & M_{\alpha}                   &           0          & M_{q}
	\end{bmatrix} \\
	\mathbf{B_{lon}} & =
	\begin{bmatrix}
	X_{\delta_{e}}  \\
	\frac{Z_{\delta_{e}}}{V_{T_{e}}} \\
	0               \\
	M_{\delta_{e}}  
	\end{bmatrix} 
	\end{align}
\end{subequations}
where the non-zero elements are known as \textit{dimensional aerodynamic derivatives} while $\theta_{e}$ is the steady-state pitch angle. The dimensional derivatives can be converted into dimensionless derivatives as shown in \eqref{eq:AdimensionaAerodynamicCoefficients} via the geometrical configuration of the aircraft, for details see \cite{stevens2015aircraft,mulder2000flight}. 
The longitudinal dynamics can be further decoupled into the \textit{Phugoid} and \textit{Short-period} mode. The Phugoid mode is normally rather slow, slightly damped, and dominates the response in $\VT$ and $\theta$, while the Short-period mode is typically fast, moderately damped, and dominates the response in $\alpha$ and $q$. For control applications, accurate knowledge of the Phugoid mode is not crucial due to the low frequency of oscillation which is compensated via feedback control, whereas the Short-period mode is crucial for stability and performance characteristics \cite{cook2007flight}. 

The lateral dynamics are described analogously by the states $\mathbf{x_{lat}} = \left[\beta~\phi~p~r\right]^\top$, which correspond to equations \eqref{eq:beta}, \eqref{eq:phi}, \eqref{eq:p}and \eqref{eq:r}. Force $Y$ and moments $L$ and $N$ are described by linear functions of $\beta,p,r$ and inputs $\mathbf{u_{lat}} = \left[ \delta_{a}~\delta_{r}\right]^\top$. The resulting matrices are given by
\begin{subequations}\label{eq:latDyn}
	\begin{align}
	\mathbf{A_{lat}} & = 
	\begin{bmatrix}
	\frac{Y_{\beta}}{V_{\mathrm{T_{e}}}} & g_{\mathrm{D}}\cos\theta_{e} & Y_{p}  & Y_{r} - V_{\mathrm{T_{e}}}  \\
	0                           &      0              &   1    & \tan \theta_{e}    \\
	L_{\beta}'                  &      0              & L_{p}' & L_{r}'             \\
	N_{\beta}'                  &      0              & N_{p}' & N_{r}' 
	\end{bmatrix} \\
	\mathbf{B_{lat}} & =
	\begin{bmatrix}
	\frac{Y_{\delta_{a}}}{V_{T_{e}}}  & \frac{Y_{\delta_{r}}}{V_{T_{e}}} \\
	0                                 &                     0            \\
	L_{\delta_{a}}'                   &              L_{\delta_{r}}'     \\
	N_{\delta_{a}}'                   &              N_{\delta_{r}}'
	\end{bmatrix}, 
	\end{align}
\end{subequations}
and their derivatives are defined in \cite{mcruer2014aircraft}. Unlike the longitudinal dynamics, the lateral motion cannot be decoupled into independent modes. They are governed by a slow \textit{Spiral} mode a fast lightly damped \textit{Dutch roll} mode, and an even faster \textit{Roll Subsidence} mode (for details see \cite{stevens2015aircraft}).

\subsection{Design of conventional maneuvers}\label{sec:DesignManeuvers}
\begin{figure}[tbhp]
	\centering
	\includegraphics[scale=0.6]{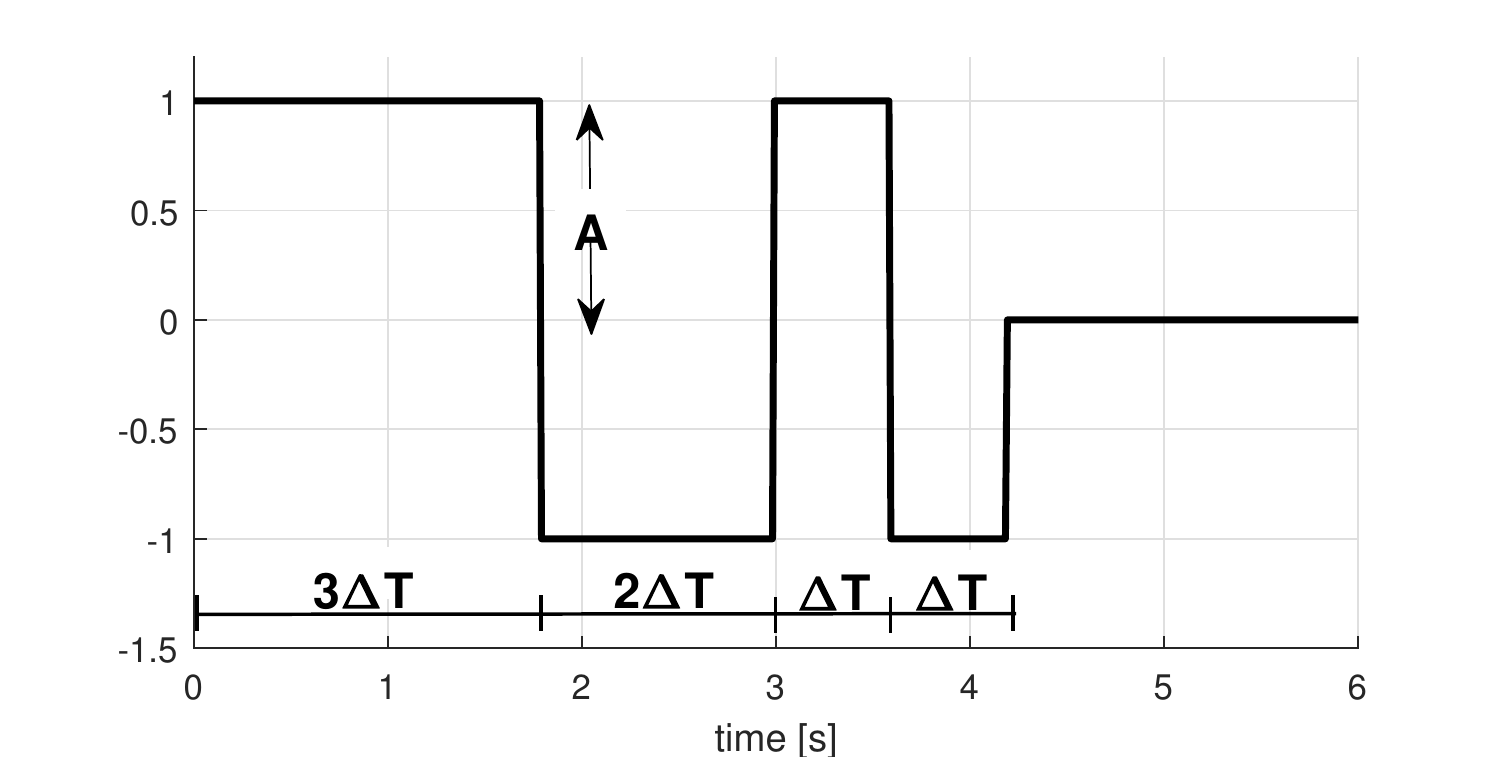} 
	\caption{Example of 3-2-1-1 maneuver with $A = 1$ and $\Delta \mathrm{T} = 0.6$.}
	\label{fig:3211maneuvre}
\end{figure}

A type of signal input for this application which is widely used in the aerospace field due to its easy implementation and good estimation performance comes from an optimization procedure of a sequence of step functions, developed by Koehler \cite{koehler1977auslegung}. The input signal has a bang-bang behavior with a duration $7\Delta \mathrm{T}$ with switching times at $t = 3\Delta \mathrm{T}$, $t = 5\Delta \mathrm{T}$, and $t = 6\Delta \mathrm{T}$ and amplitude $A$. For this reason, such an input signal is called a \textit{3-2-1-1 maneuver} (see figure~\ref{fig:3211maneuvre}).

In \cite{mulder1994identification}, it was shown that the 3-2-1-1 maneuver provides the best estimation accuracy for both aircraft longitudinal and lateral dynamics among \textit{Doublets}, \textit{Mehra}, \textit{Schulz} and \textit{DUT} input signals. Yet, only Doublets and 3-2-1-1 input signals provide sufficient system excitation for identification of system responses with frequencies above $1\,\mathrm{Hz}$, though the 3-2-1-1 maneuver embraces much higher frequencies compared to Doublets.  
Finally, 3-2-1-1 maneuvers can be chosen through both a qualitative consideration in the frequency domain \cite{marchand1977untersuchung} and a \textit{trial-and-error} approach in order to ensure that the system response is within the flight envelope.

\subsection{Design of optimized maneuvers}\label{sec:DesignOptimizedManeuvers}
Another type of signal input implemented in this work is obtained by solving a time domain model-based \ac{OED} problem that aims to obtain more accurate parameter estimates while enforcing safety constraints \cite{licitra2017oed}.

The main idea of \ac{OED} is to use as objective of an optimization problem a function $\mathbf{\Psi}(\cdot)$ of the Fisher information matrix $\mathbf{F}$ which is given by  
\begin{equation}
\mathbf{F} = \sum_{i = 1}^{N} \left[ \left( \frac{\partial \mathbf{y}(i)}{\partial \mathbf{\tilde{p}}} \right)^{T} \mathbf{\Sigma_{y}}^{-1} \left( \frac{\partial \mathbf{y}(i)}{\partial \mathbf{\tilde{p}}} \right)\right].
\label{eq:FisherMatrix}
\end{equation}
with $\mathbf{y} \in \mathbb{R}^{n_{\mathbf{y}}}$ the output states sampled in $N$ measurements, a priori parameters $\mathbf{\tilde{p}} \in \mathbb{R}^{n_{\mathbf{p}}}$ and $\mathbf{\Sigma_{y}} \in \mathbb{R}^{n_\mathbf{y} \times n_\mathbf{y}}$ the measurements noise covariance matrix.
A general model-based OED problem which considers input $\mathbf{u}(t)$, state constraints $\mathbf{x}(t)$, time length $\mathrm{T}$ and subject to a mathematical model expressed as an \ac{ODE} can be formulated as
\begin{subequations}\label{eq:DOE_TC}
	\begin{align}
	\underset{\mathbf{x}(\cdot),\mathbf{u}(\cdot)}{\text{minimize}} \hspace{1cm} & \mathbf{\Psi} \left(\mathbf{F} \left[\mathbf{x}(\cdot), \mathbf{u}(\cdot), \mathbf{\tilde{p}}\right] \right) \label{eq:OCP_objective} \\ 
	\text{subject to:} \hspace{1cm} & \dot{\mathbf{x}}(t) = \mathbf{f} \left(\mathbf{x}(t), \mathbf{u}(t), \mathbf{\tilde{p}}\right), t \in \left [ 0,\mathrm{T}  \right ]\label{eq:OCP_ode}  \\ 
	& \mathbf{x}(0) = \mathbf{x}_{0}, \\
	& \mathbf{u}_{\mathrm{min}} \leq \mathbf{u}(t) \leq \mathbf{u}_{\mathrm{max}}  \hspace{0.49cm}, t \in  \left [ 0,\mathrm{T}  \right ] \label{eq:OCP_control_bounds} \\
	& \mathbf{x}_{\mathrm{min}} \leq \mathbf{x}(t)  \leq \mathbf{x}_{\mathrm{max}} \hspace{0.51cm},  t \in \left [ 0,\mathrm{T} \right] \label{eq:OCP_state_bounds}
	\end{align}
\end{subequations}
For further details refer to \cite{licitra2017optinput,OptSysIDvideo}.

\subsection{Baseline model}\label{sec:baselineModel}
Both 3-2-1-1 and \ac{OED}-based maneuvers need to be designed using a baseline (a priori) model with reasonable accuracy so as to both have an first insight about the estimation performance and ensure that the system response evolves within the flight envelope. 

Various methods can be applied to obtain a priori models. If the airframe is similar to an existing aircraft, its model can be scaled. For instance, the Digital DATCOM \cite{hoak1975usaf} is a purely empirical guide to estimating aerodynamic derivatives based on aircraft configuration and the experience of engineers. If the airfoils and aircraft configuration are new, one can perform analysis via the lifting line method \cite{versteeg2007introduction}, \ac{CFD} \cite{anderson2017fundamentals}, wind-tunnel tests or previous flight tests. Depending on the available resources, combinations of these methods can be used. In this work, a priori models are retrieved from lifting line method \cite{licitra2017pe}. 

A steady wing-level flight condition is considered with trimmed airspeed $\VTe = 20\,\mathrm{m/s}$. The equilibrium point is held for $\delta_{e} = -1.5 \, ^\circ, \alpha_{\mathrm{e}} = - 0.4 \, ^\circ$ and $\theta_{\mathrm{e}} = - 4.5 \, ^\circ$ with the other states equal to zero. Subsequently the system is linearized and the longitudinal dynamics are taken into account as in \eqref{eq:longDyn} with dimensional a priori derivatives shown in Table~\ref{tab:LongitudinalDerivatives}.
\begin{table}[tbhp]
	\caption{A priori longitudinal dimensional aerodynamic derivatives}
	\label{tab:LongitudinalDerivatives}
	\centering
	\begin{tabular}{ lr lr lr }\hline
		\bf X-axis       &\bf Value  & \bf Z-axis              &\bf Value &\bf M-axis        &\bf Value \\\hline
		$X_{V}$          &  -0.147   & $Z_{V}$                 & -0.060   & $M_{V}$          &   0.0    \\
		$X_{\alpha}$     &   7.920   & $Z_{\alpha} / \VTe$     & -4.400   & $M_{\alpha}$     &  -6.180  \\
		$X_{q}$          &  -0.163   & $Z_{q}$                 &  0.896   & $M_{q}$          &  -1.767  \\
		$X_{\delta_{e}}$ &  -0.232   & $Z_{\delta_{e}} / \VTe$ & -0.283   & $M_{\delta_{e}}$ & -10.668  \\\hline		
	\end{tabular}
\end{table}

As mentioned in section~\ref{sec:steadyConditionAndDecoupling}, the a priori models provide an insight into the general characteristics of the aircraft behavior via modal analysis \cite{stevens2015aircraft}. In Table~\ref{tab:AprioriModes} the a priori aircraft modes relative to the longitudinal dynamics are provided in terms of natural frequencies $\omega_{n}$, damping ratios $\delta$, constant times $\tau$, overshoots in percentage $S_{\%}$ and period of oscillations $P_{\mathrm{O}}$.

\begin{figure}[tbhp]
	\centering
	\includegraphics[scale=0.60]{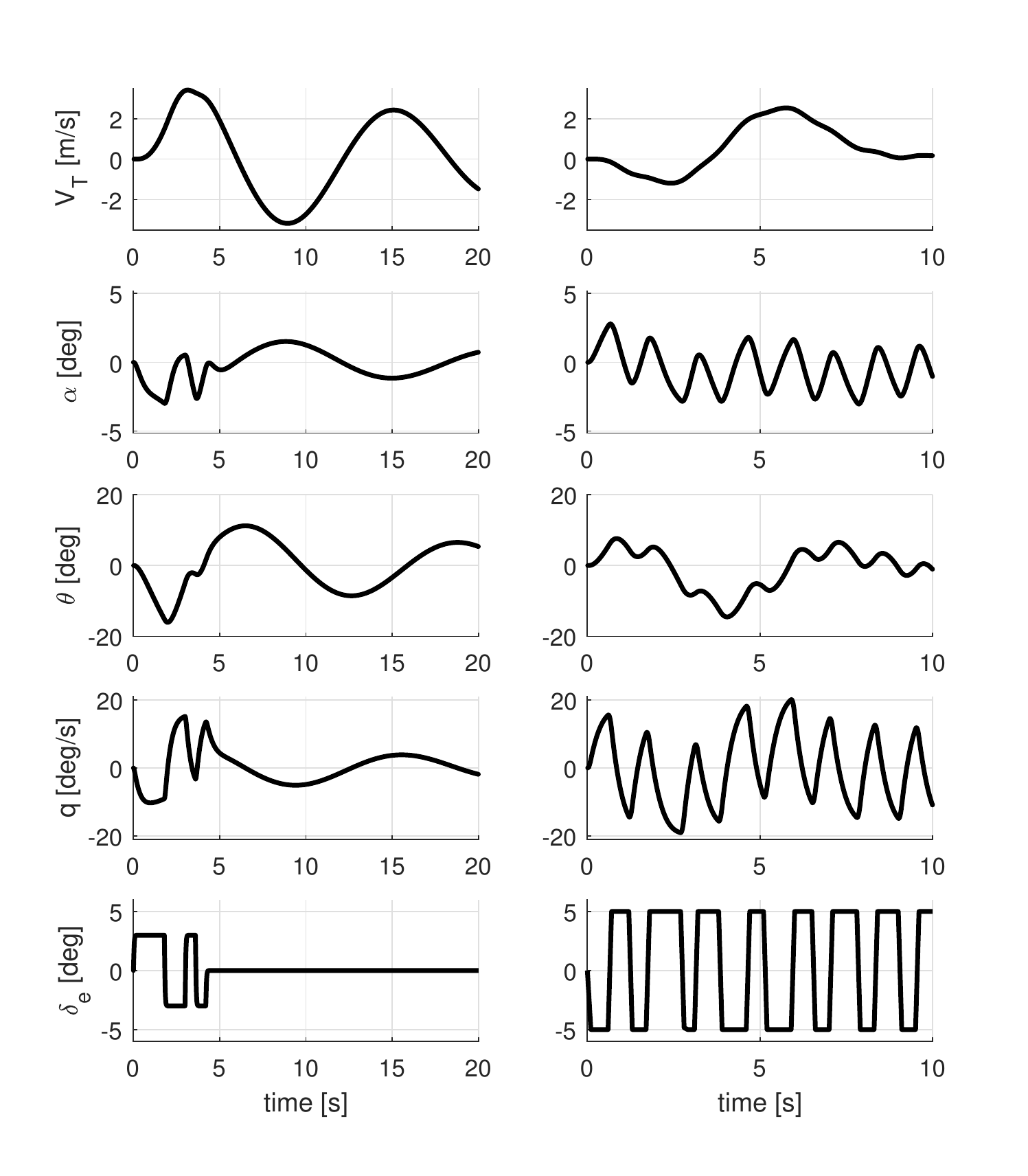} 
	\caption{A priori system response with both 3-2-1-1 and optimized maneuvers.}
	\label{fig:AprioriSimulation}
\end{figure}

The modal analysis suggests to design experiments with time duration longer than $12.067 \; \mathrm{s}$ so as to provide sufficient excitation in the frequency range where the expected Phugoid mode is defined. 
Figure~\ref{fig:AprioriSimulation} shows the candidate maneuvers with the corresponding a priori system response. Note that the experiment relative to the 3-2-1-1 maneuver has a time length of $20 \; \mathrm{s} $ whereas the optimized experiments length are set to $10\; \mathrm{s}$ to ensure the full sequence is completed in the available flight test area taking into account variations in the wind conditions on the flight test day(s).

%the optimized experiment is set to $10$[\rm{s}] due to limitation of the flight test field for that specific flight campaign.
%since it is likely that such an \textit{aggressive} maneuver performed in open-loop might lead the aircraft to unsafe regions.

\begin{table}[tbhp]
	\caption{A priori longitudinal modes}
	\label{tab:AprioriModes}
	\centering
	\begin{tabular}{cccc}
		\hline       
		\bf Mode         & \bf Short-period & \bf Phugoid&          \\\hline
		$\omega_{n}$     &   3.939      &  0.521 & $\mathrm{rad/s}$ \\
		$\tau$           &   0.254      &  1.920 & $\mathrm{s}$     \\
		$\delta$         &   0.789      &  0.031 & $ - $            \\
		$S_{\%}$         &   1.768      & 90.831 & \%               \\
		$P_{\mathrm{O}}$ &   2.596      & 12.067 & $\mathrm{s}$     \\\hline
	\end{tabular}
\end{table}

%For this reason, \textit{flight envelope limit detection} algorithms should be programmed in the \ac{FCC} in order to avoid flight envelope violation e.\,g., due to significant inaccuracies of the a priori models or unexpected gust occurring during the open loop-phase

Historically, aircraft system identification has been performed using a pilot to provide input sequences.  In this work, the input sequences are performed autonomously. The flight control computer monitors the aircraft response and aborts the maneuver whether the predetermined flight envelope boundaries are violated \cite{licitra2017oed,licitra2017pe,licitra2017optinput}.  

\subsection{Preliminary analysis}\label{sec:PreliminaryAnalysis}
One way to assess the estimation accuracy that a given maneuver can provide is by the \ac{CRLB}, i.e., the theoretical lower limits for parameter standard errors $\sigma$ using an efficient and asymptotically unbiased estimator, such as maximum likelihood \cite{morelli1993oed}. 
A performance analysis of signal inputs computed via the \ac{CRLB} isolates the merits of the input design from the merits of the parameter estimation algorithm used to extract the aerodynamic derivatives from the flight data \cite{klein1990optimal}. 
The \ac{CRLB} depends on the diagonal entries of the Fisher information matrix $\mathbf{F}$ \eqref{eq:FisherMatrix} which is formally \cite{soderstrom2001identification}
\begin{equation}\label{eq:std-CR}
\sigma_{i} \geq \mathrm{CRLB}_{i} = \frac{1}{\sqrt{\mathbf{F}_{ii}}}.
\end{equation}
Experience has shown that a factor of 2 can be used so as to obtain an approximation of the parameter standard error \cite{tischler2006aircraft}, resulting in  
\begin{equation}\label{eq:2std-2CR}
\sigma_{i} \approx 2 \cdot \mathrm{CRLB}_{i} = \frac{2}{\sqrt{\mathbf{F}_{ii}}}.
\end{equation}
Finally, note that the inverse of the Fisher information matrix $\mathbf{F}^{-1}$ corresponds to the covariance matrix of the estimated parameters $\mathbf{\Sigma_{\mathbf{p}}} \in \mathbb{R}^{n_{\mathbf{p}} \times n_{\mathbf{p}}}$. Table~\ref{tab:CRlongDyn} gathers the 2\ac{CRLB} values in percentage for the system responses shown in figure~\ref{fig:AprioriSimulation} and using the sensors noise standard deviation $\sigma_{\mathbf{y}}$ collected in table~\ref{tab:std_sensors}. 

\begin{table}[tbhp]
	\caption{Dimensional aerodynamic longitudinal derivatives with corresponding expected estimation accuracy via $\mathrm{2CRLB}$.}
	\label{tab:CRlongDyn}
	\centering
	\begin{tabular}{lrr}\hline
		\bf Derivatives       &\bf Value & $\mathbf{2CRLB_{\%}}$\\\hline
		$X_{V}$               &  -0.064  &   25.27   \\
		$X_{\alpha}$          &   8.635  &   34.88   \\
		$X_{q}$               &  -0.153  &  336.67   \\			
		$X_{\delta_{e}}$      &  -0.173  &  291.62   \\			
		$Z_{V}$               &  -0.050  &    1.07   \\
		$Z_{\alpha}/\VTe$     &  -4.222  &    1.06   \\
		$Z_{q}$               &   0.897  &    1.08   \\			                    
		$Z_{\delta_{e}}/\VTe$ &  -0.340  &    4.18   \\
		$M_{\alpha}$          &  -7.671  &    0.14   \\
		$M_{q}$               &  -1.963  &    0.10   \\	   
		$M_{\delta_{e}}$      & -17.939  &    0.02   \\\hline		
	\end{tabular}
\end{table}
\begin{table}[tbhp]
	\caption{Sensors noise standard deviation $\sigma_{y}$}
	\label{tab:std_sensors}
	\centering
	\begin{tabular}{lccc}\hline
		\bf Sensor           & \bf Variable         & $\mathbf{\sigma_{y}}$ & \bf Unit             \\ \hline
		five hole pitot tube & $V_{T}$              &    $1.0$     & $\mathrm{m/s}$  \\
		five hole pitot tube & $(\alpha,\beta)$     &    $0.5$     & $\mathrm{deg}$   \\
		IMU	                 & $(\phi,\theta,\psi)$ &    $0.1$     & $\mathrm{deg}$   \\ 
		IMU	                 & $(p,q,r)$            &    $0.1$     & $\mathrm{deg/s}$ \\\hline		
	\end{tabular}
\end{table}
The results indicate that the dimensional aerodynamic derivatives relative to the Phugoid mode, i.\,e., $X_{q}$, $X_{\delta_{e}}$ which correspond to the dimensionless one 
$C_{X_{q}}$, $C_{X_{\delta_{e}}}$ are subject to high uncertainty. 
High values of \ac{CRLB} indicate that either the $i^\mathrm{th}$ parameter is physically insignificant with respect to the measured aircraft response or there is a correlation between parameters, i.\,e., these parameters can vary together, making their individual values difficult to determine \cite{soderstrom2001identification}. In this case, it turns out that $X_{V}$ provides a negligible contribution to the aircraft response whereas correlation occurs between $X_{\alpha}$, $X_{q}$ and $X_{\delta_{e}}$. 

To overcome this issue, one might fix the parameters associated to the Phugoid mode with their a priori values though, errors in the form of a low-frequency model mismatch could arise in the identified model. Nevertheless, accurate knowledge of the Phugoid mode is not crucial to due its slow motion, which can be easily handled by a pilot or a control system \cite{dorobantu2013system}. On the other hand, high estimation accuracy is required for the Short-period mode which is given by $Z_{V}$, $Z_{\alpha}$, $Z_{q}$, $Z_{\delta_{e}}$, $M_{\alpha}$, $M_{q}$ and $M_{\delta_{e}}$ since longitudinal stability and performance characteristics primarily depend on the accuracy of the Short-period mode \cite{cook2007flight}. 

\section{Experimental Data}\label{sec:ExperimentalData}
\begin{figure}[tbhp]
	\centering
	\includegraphics[scale=0.6]{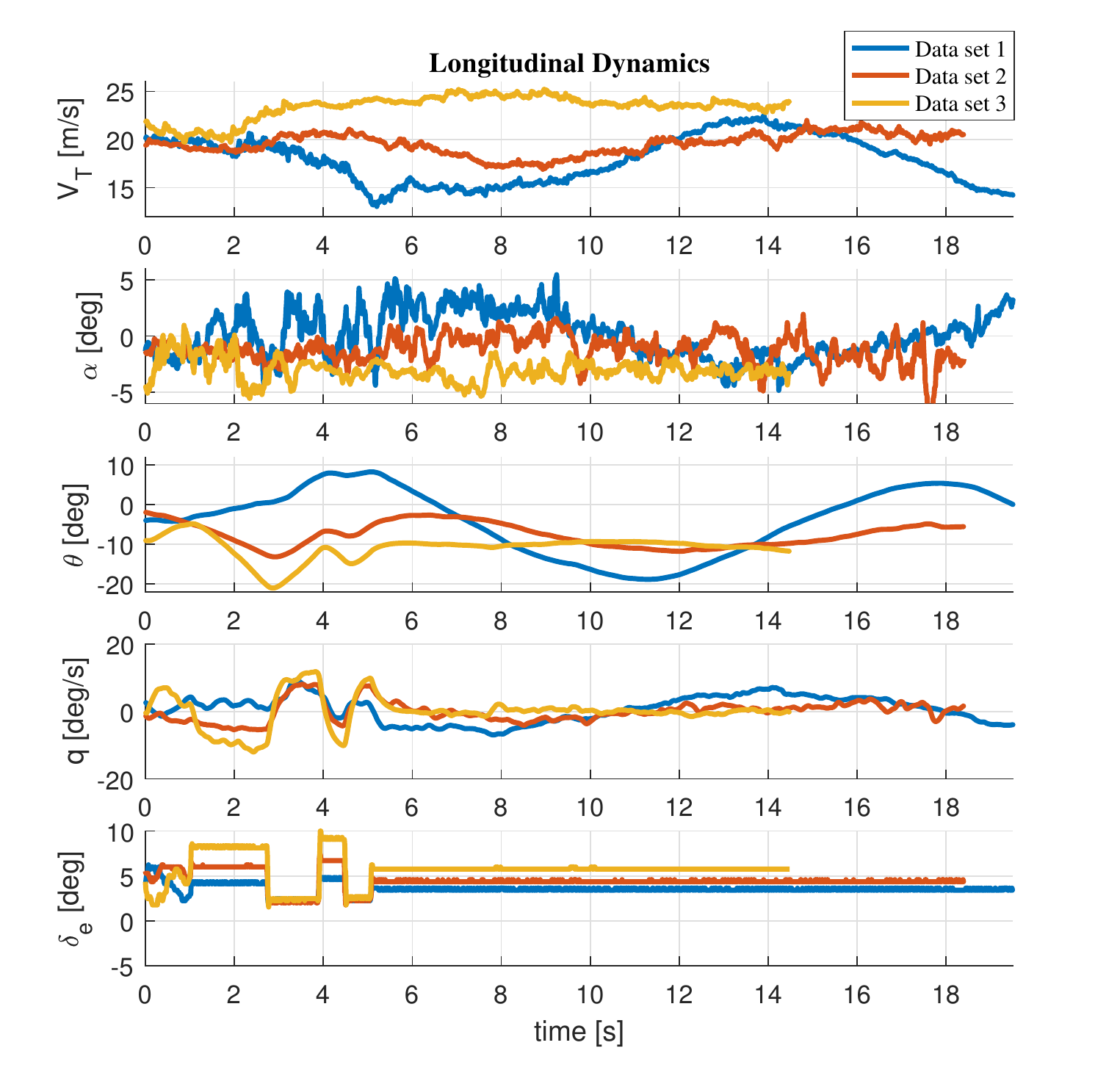} 
	\caption{Three experimental data sets obtained through conventional maneuvers. Average wind speed $\approx 7\,\mathrm{m/s}$.}
	\label{fig:EstimationDataSet_LongDyn_3211}
\end{figure}

\begin{figure}[tbhp]
	\centering
	\includegraphics[scale=0.6]{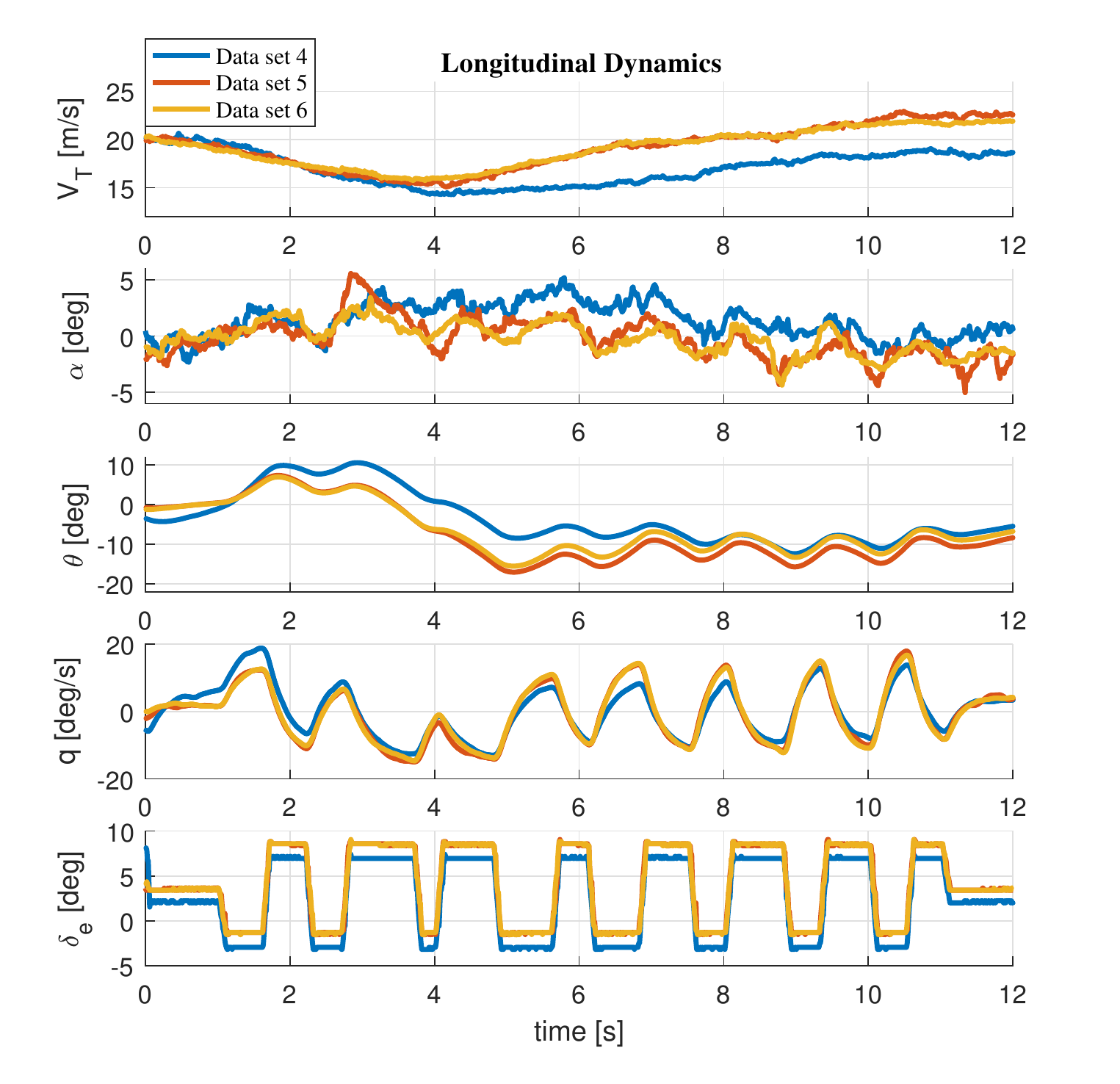} 
	\caption{Three experimental data set obtained through optimized maneuvers. Average wind speed $\approx 2\,\mathrm{m/s}$.}
	\label{fig:EstimationDataSet_LongDyn_Opt}
\end{figure}

Six experimental data are collected from two different flight test campaign. Three experiments are performed with conventional maneuvers 3-2-1-1 shown in figure~\ref{fig:EstimationDataSet_LongDyn_3211} with an average (estimated) wind speed $\approx 7 \, \mathrm{m/s}$ whereas other three experiments are collected using the \ac{OED}-based maneuvers and shown in figure~\ref{fig:EstimationDataSet_LongDyn_Opt} with average wind speed $\approx 2\,\mathrm{m/s}$. 

In figure~\ref{fig:EstimationDataSet_LongDyn_Opt} one can observe the decoupling between the Phugoid mode which dominates the airspeed $\VT$ and pitch $\theta$ responses, with the fast changes on the angle of attack $\alpha$ and pitch rate $q$ coming from the Short-period mode.  
Comparing figure~\ref{fig:EstimationDataSet_LongDyn_3211} with figure~\ref{fig:EstimationDataSet_LongDyn_Opt}, it is possible to discern the turbulence effect on the angle of attack $\alpha$ and pitch rate $q$ response.  
This is not surprising since turbulences increase consistently with the wind speed.

\begin{figure}[tbhp]
	\centering
	\includegraphics[scale=0.60]{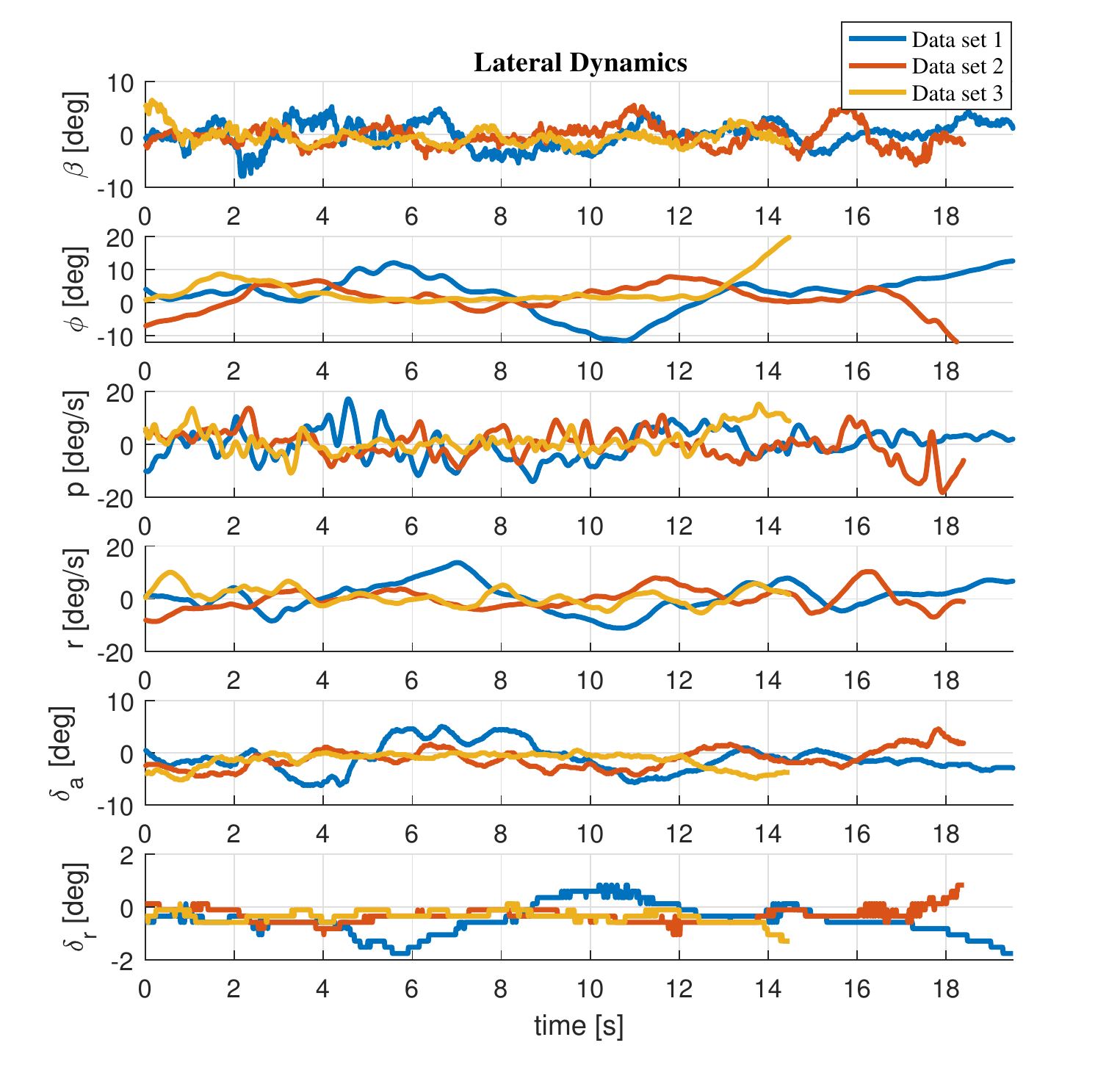} 
	\caption{Stabilization of lateral dynamics by $\delta_{a}$ and $\delta_{r}$ during excitation signal along the longitudinal dynamics via conventional maneuvers. Average wind speed $\approx 7 \,\mathrm{m/s}$.}
	\label{fig:EstimationDataSet_LatDyn_3211}
\end{figure}

\begin{figure}[tbhp]
	\centering
	\includegraphics[scale=0.60]{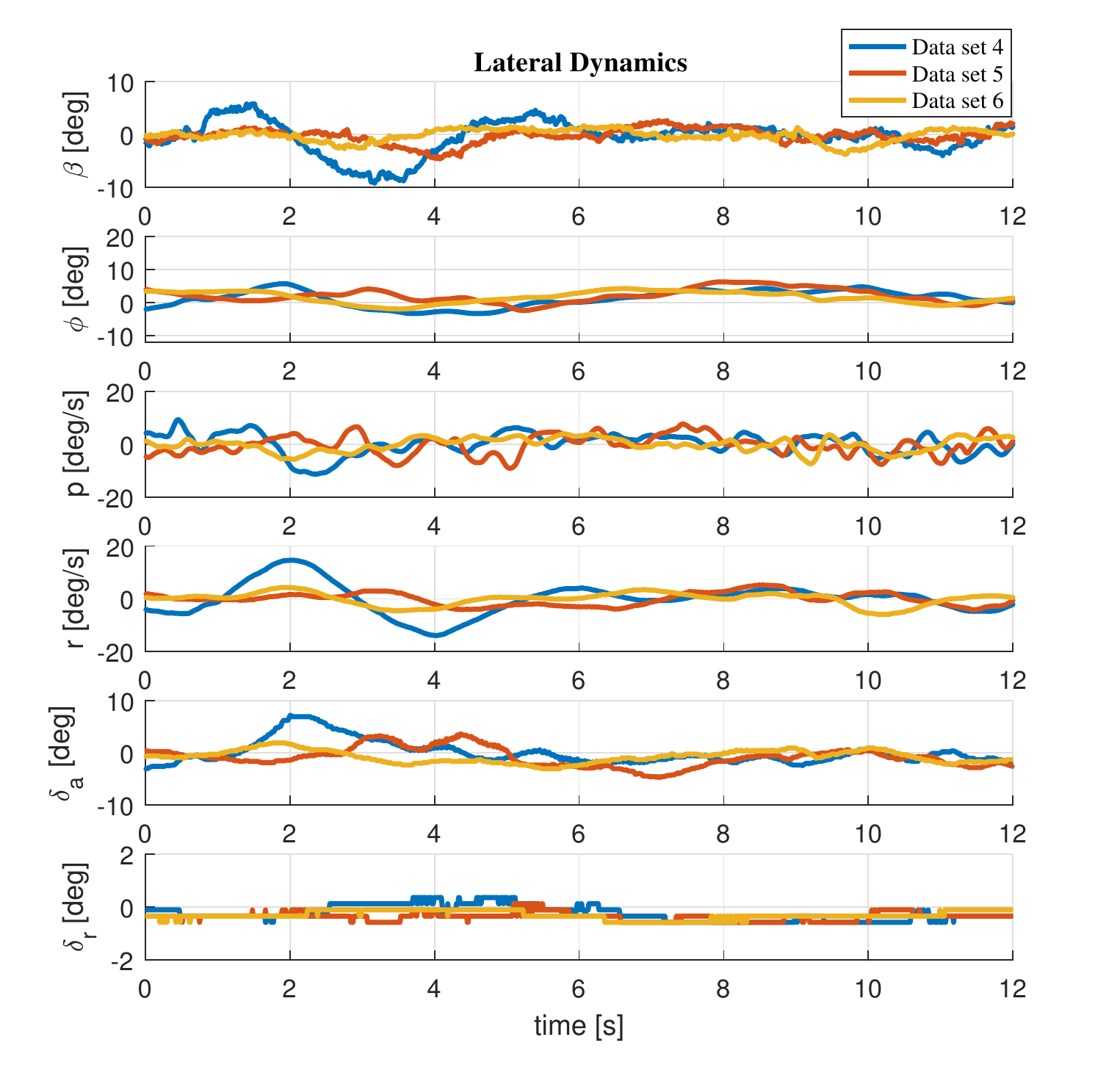} 
	\caption{Stabilization of lateral dynamics by $\delta_{a}$ and $\delta_{r}$ during excitation signal along the longitudinal dynamics via optimized maneuvers. Average wind speed $\approx 2 \,\mathrm{m/s}$.}
	\label{fig:EstimationDataSet_LatDyn_Opt}
\end{figure}

As mentioned in section \ref{sec:steadyConditionAndDecoupling}, during the excitation of the longitudinal dynamics, the lateral motion need to be stabilized throughout the entire  experiment via aileron $\delta_{a}$ and rudder $\delta_{r}$ deflection. Figures~\ref{fig:EstimationDataSet_LatDyn_3211} and \ref{fig:EstimationDataSet_LatDyn_Opt} show the lateral dynamics relative to the conventional and optimized experiments, respectively. Also in this case, it is clearly shown how the turbulence effect acts on the aircraft dynamics.\\ 
More precisely, in figure~\ref{fig:EstimationDataSet_LatDyn_3211} one can observe how the roll rate $p$ and roll angle $\phi$ appear sensitive to the turbulence which involve a major control effort from the aileron deflection $\delta_{a}$ in order to both stabilize this axis and prevent flight envelope violation. 


\section{Formulation of multiple experiment model based parameter estimation}\label{sec:MBPEformulation}
Whenever parameter estimation is intended for identification of aircraft dynamics, multiple experiments are usually required to deal with the following issues \cite{morelli2006practical}
\begin{itemize}
	\item Reduce the effects of sensor biases as well as atmospheric turbulence on estimation results;
	\item individual maneuvers usually have good information content only for a subset of parameters, while multiple maneuvers combined can provide better information w.r.t the complete set of parameters;
	\item the flight test area and operating safety case restricts the flight paths that can be flown, limiting the available duration of any particular maneuver.
\end{itemize}
A standard approach is to retrieve the estimated parameters via data fitting for each independent experiment and subsequently weight them w.r.t. to their inverse (estimated) parameter covariance matrix $\mathbf{\Sigma_{\mathbf{p}}}$ \cite{ljung1998system}. However, such method might lead to wrong results whenever computed $\mathbf{\Sigma_{\mathbf{p}}}$ are not reliable \cite{licitra2017pe}. 

Furthermore, in equation~\eqref{eq:MathematicalModel} one can observe that angular acceleration measurements as well as rate of changes in the airspeed, Euler and aerodynamic angles are required in order to estimate aerodynamic properties. Usually, these quantities are not measured though, they can be retrieved by numerical differentiation methods, which are rather noisy \cite{morelli2006practical}. Consequently, signal distortion may arise degrading the overall estimation performance. Within this scenario, multiple experiments \ac{MBPE} algorithms appear a reasonable choice for estimation of aerodynamic derivatives. 

In this context, let us consider a mathematical model defined as a set of \ac{ODEs}
\begin{subequations}\label{eq:ODE}
	\begin{align}	
	\dot{\mathbf{x}}(t) & = \mathbf{f}(\mathbf{x}(t), \mathbf{u}(t), \mathbf{p},t) \label{eq:f} \\
	\mathbf{y}(t)       & = \mathbf{h}(\mathbf{x}(t), \mathbf{u}(t), \mathbf{p},t) + \mathbf{\epsilon}(t) \label{eq:h}
	\end{align}
\end{subequations}
with differential states $\mathbf{x} \in \mathbb{R}^{n_{\mathbf{x}}}$, output state $\mathbf{y} \in \mathbb{R}^{n_{\mathbf{y}}}$ noise-free control inputs $\mathbf{u} \in \mathbb{R}^{n_{\mathbf{u}}}$, parameters $\mathbf{p} \in \mathbb{R}^{n_{\mathbf{p}}}$, and time $t$. The measurement values $\mathbf{y}$ are polluted by additive, zero-mean Gaussian noise $\eta \left(0,\mathbf{\Sigma_{y}}\right)$ with $\mathbf{\Sigma_{y}}$ the covariance of noise measurements.  

A multiple experiments \ac{MBPE} problem can be first stated using an \ac{OCP} perspective in continuous time as follows \cite{licitra2017pe}
% Parameter Estimation formulation in Continuous Time Domain
\begin{subequations}\label{eq:OPC_TC}
	\begin{align}
	& \underset{\mathbf{p\left(\cdot\right)}}{\text{minimize}} & & \sum_{i=1}^{N_{\mathrm{e}}}  \int_{0}^{T^{i}}  \left \|   \hat{\mathbf{y}}^{i}(t) - \mathbf{h}\left ( \mathbf{x}^{i}(t),\hat{\mathbf{u}}^{i}(t),\mathbf{p} \right ) \right \|^{2}_{\mathbf{\Sigma_{y}}^{-1}} \mathrm{d}t \\
	& \text{subject to}
	& & \dot{\mathbf{x}}^{i}(t)= \mathbf{f}(\mathbf{x}^{i}(t),\hat{\mathbf{u}}^{i}(t),\mathbf{p},t)\\
	& & & t \in \left [ 0,T^{i} \right ],\; i \in \mathbb{Z}_{1}^{N_{\mathrm{e}}}
	\end{align}
\end{subequations}
with $N_{\mathrm{e}}$ number of experiments, $\hat{\mathbf{u}}^{i}(t)$ and $\hat{\mathbf{y}}^{i}(t)$ respectively input, output measurements for $i^{\mathrm{th}}$ experiment running for a duration $T^{i}$. Using \textit{direct methods} \cite{diehl2014ocp}, the optimization problem \eqref{eq:OPC_TC} can be transformed into a finite dimensional \ac{NLP} which can then be solved by numerical optimization methods.  
In this work, a \textit{direct multiple shooting} approach is chosen due to its stability w.r.t. the initial guess compared to a \textit{single shooting} strategy \cite{bock1984multiple}. 

In order to implement a multiple shooting algorithm, let us define an equidistant grid over the experiment consisting in the collection of time points $t_{k}$, where $t_{k+1}-t_{k}=\frac{T^{i}}{N^{i}_{\mathrm{m}}}:=T_{s},\, \forall i=0,\ldots,N_{\mathrm{e}}$ with $N_{\mathrm{m}}^{i}$ the number of measurements for the $i^{\mathrm{th}}$ data set, assuming implicitly that the measurements are collected with a fixed sample time $T_{\mathrm{s}}$. Additionally, we consider a piecewise constant control parametrization $u(\tau)=u_{k}$ for $\tau \in [t_{k},t_{k+1})$. A function $\mathbf{\Pi}(.)$ over each shooting interval is given, which represents a numerical approximation for the solution $x_{k+1}$ of the following \ac{IVP}

\begin{equation}\label{eq:IVP}
\begin{aligned}
\dot{\mathbf{x}}(\tau)= \mathbf{f}(\mathbf{x}(\tau), \mathbf{u}_{k}, \mathbf{p},\tau), \; \tau \in [t_{k},t_{k+1}].
\end{aligned}
\end{equation}  
Such function is evaluated numerically via integration methods, such as the Runge-Kutta of order 4 (RK4) as implemented in this work. Therefore, the \ac{OCP} \eqref{eq:OPC_TC} can be translated into the \ac{NLP}

% Parameter Estimation formulation in Discrete Time Domain [Multiple Shooting]
\begin{subequations}\label{eq:OPC_TD_MS}
	\begin{align}
	& \underset{\mathbf{p},\mathbf{X}}{\text{minimize}} & & \sum_{i=1}^{N_{\mathrm{e}}}  \sum_{k=0}^{N_{\mathrm{m}}}  \left\| \hat{\mathbf{y}}^{i}_{k} - \mathbf{h}\left (\mathbf{x}^{i}_{k},\hat{\mathbf{u}}^{i}_{k},\mathbf{p} \right ) \right\|^{2}_{\mathbf{\Sigma_{y}}^{-1}}  \label{eq:cost function} \\
	& \text{subject to}
	& & \mathbf{x}^{i}_{k+1} - \mathbf{\Pi}(\mathbf{x}^{i}_{k},\hat{\mathbf{u}}^{i}_{k},\mathbf{p}) = 0 \label{eq:continuity condition}\\
	& & & k = 0, 1, ..., N_{\mathrm{m}} - 1, \;\;  i \in \mathbb{Z}_{1}^{N_{\mathrm{e}}}
	\end{align}
\end{subequations}
where $\mathbf{X} \in \mathbb{R}^{n_{\mathbf{X}}}$ with $n_{\mathbf{X}} = \sum_{i=1}^{N_{\mathrm{e}}} n_{x} (N_{\mathrm{m}}^{i}-1)$ and sorted as
\begin{equation}\label{eq:X}
\begin{aligned}
\mathbf{X} = [x^{1}_{0}, \ldots , x^{1}_{N_{\mathrm{m}}^{1}-1}, \ldots , x^{N_{\mathrm{e}}}_{0} , \ldots , x^{N_{\mathrm{e}}}_{N_{\mathrm{m}}^{\mathrm{e}}-1}]^{T}
\end{aligned}
\end{equation}
in order to create a block diagonal structure on the \ac{NLP} formulation and especially in the equality constraints \eqref{eq:continuity condition}. Notice that in \eqref{eq:X} the number of measurements $N_{m}$ are assumed different for each $i^{\mathrm{th}}$ experiment. 

Finally, the \ac{NLP} initialization can be chosen from, e.g., previous estimates of $\mathbf{p}$ while $\mathbf{X}$ can be initialized using the measurements $\hat{\mathbf{y}}$ and/or estimates of the state $\mathbf{x}$. For further details refer to \cite{diehl2014ocp,bock2013model}.   

\section{Parameter estimation results}
In this section, the \ac{PE} is carried out on the experimentally obtained data. The estimation results are subsequently assessed via a time domain model validation approach. 
\subsection{Data fitting}\label{sec:dataFitting}
Within this work, the multiple experiment \ac{MBPE} algorithm is implemented using \textsc{CasADi} \cite{Andersson2013b} in \textsc{Matlab} environment.
The system dynamics \eqref{eq:f} taken into account are the non-linear longitudinal motion expressed in \eqref{eq:Vt},\eqref{eq:alpha},\eqref{eq:theta},\eqref{eq:q} with differential states 
\begin{equation}
\mathbf{x}(t) = [\VT(t)~\alpha(t)~\theta(t)~q(t)]^{\top}
\end{equation}
assuming steady wing-level flight condition, i.e., $\beta = \phi = p = r = 0$. The unknown parameters are
\begin{equation}
\mathbf{p} =
\begin{pmatrix} 
C_{X_{0}} & C_{X_{\alpha}} & C_{X_{q}} & C_{X_{\delta_{e}}} \\ 
C_{Z_{0}} & C_{Z_{\alpha}} & C_{Z_{q}} & C_{Z_{\delta_{e}}} \\
C_{m_{0}} & C_{m_{\alpha}} & C_{m_{q}} & C_{m_{\delta_{e}}} 
\end{pmatrix}
\in \mathbb{R}^{3 \times 4}
\end{equation}
and control input equal to 
\begin{equation}
\mathbf{u}(t) = \delta_{e}(t).
\end{equation}
whereas the output states \eqref{eq:h} are simply given by
\begin{equation}
\mathbf{y}(t) = \mathbf{x}(t) + \mathbf{\epsilon}(t).
\end{equation}
The continuous-time optimization problem \eqref{eq:OPC_TC} is subsequently discretized and formulated as a \ac{NLP} using direct multiple shooting. The resulting \ac{NLP} is solved via \textsc{IPOPT} \cite{Waechter2006} with linear solver \textsc{MA27} \cite{HSL2017}.
Finally, the optimization problem \eqref{eq:OPC_TD_MS} is initialized using the baseline model described in section \ref{sec:baselineModel} for $\mathbf{p}$ and $\mathbf{X}$ with the real output measurements $\hat{\mathbf{y}}^{i}, i \in \mathbb{Z}_{1}^{N_{\mathrm{e}}}$.

\begin{figure}[tbhp]
	\centering
	\includegraphics[scale=0.60]{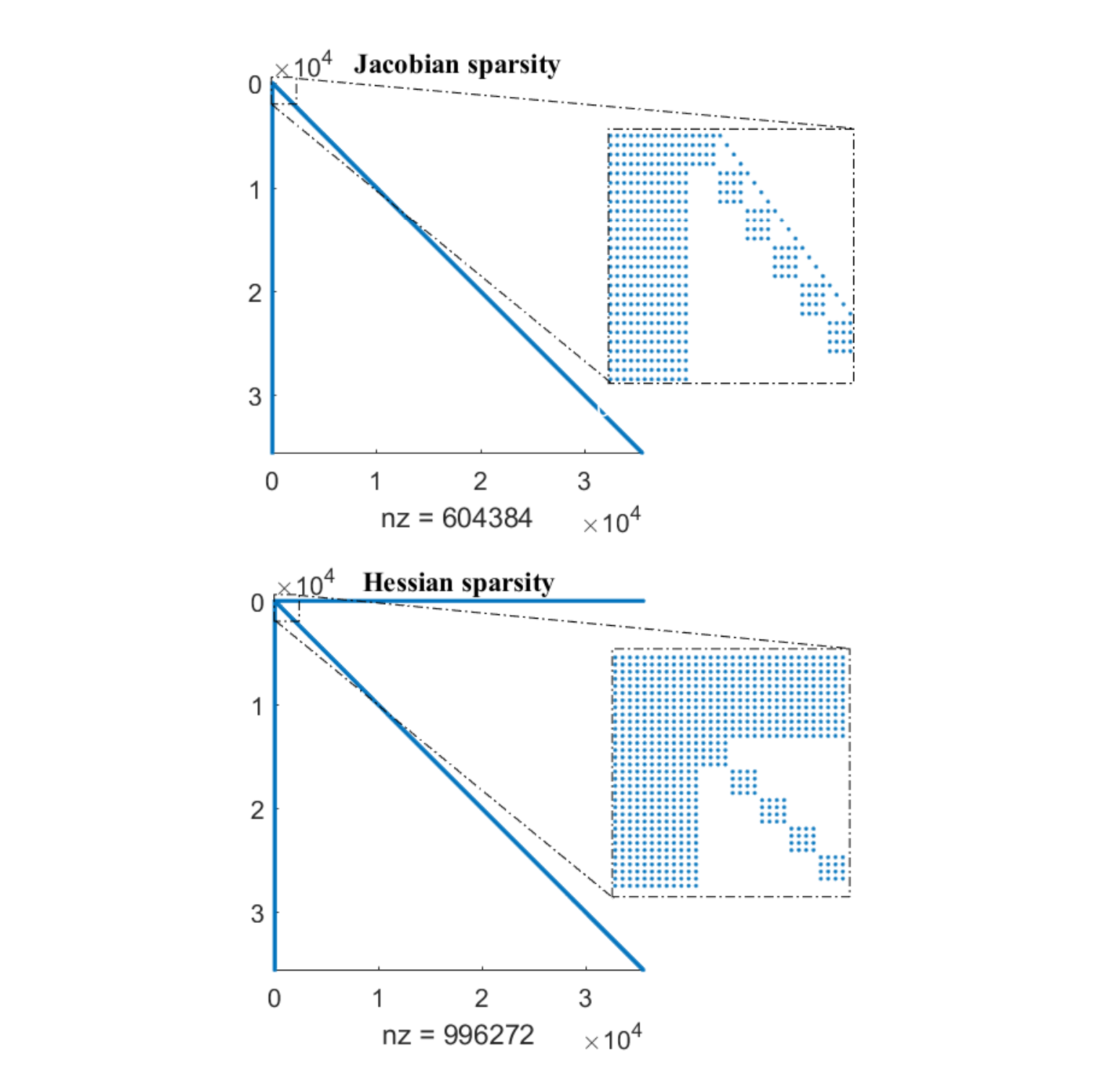} 
	\caption{Jacobian and Hessian Sparsity of the \ac{NLP}.}
	\label{fig:NLPsparsity}
\end{figure}

The data fitting is carried out simultaneously for all the experimental data set shown in section \ref{sec:ExperimentalData} with total number of optimization variables
\begin{equation}
n_{\mathbf{opt}} = n_{\mathbf{p}} + n_{\mathbf{X}} = 12 + 35564 = 35576.
\end{equation} 
\textsc{CasADi} discovers the structure and computes the full sparse Jacobian and Hessian with a minimal of algorithmic differentiation sweeps (see figure~\ref{fig:NLPsparsity}). \textsc{CasADi's} for-loop equivalents are used to efficiently build up the large number of shooting constraints \eqref{eq:continuity condition}. Furthermore, since this application requires a large number of control intervals, the \textsc{CasADi} \texttt{map} functionality was used to achieve a memory-lean computational graph. Using this proposed implementation, the \ac{NLP} is solved within 28 iterations of \textsc{IPOPT}.

\begin{figure}[tbhp]
	\centering
	\includegraphics[scale=0.6]{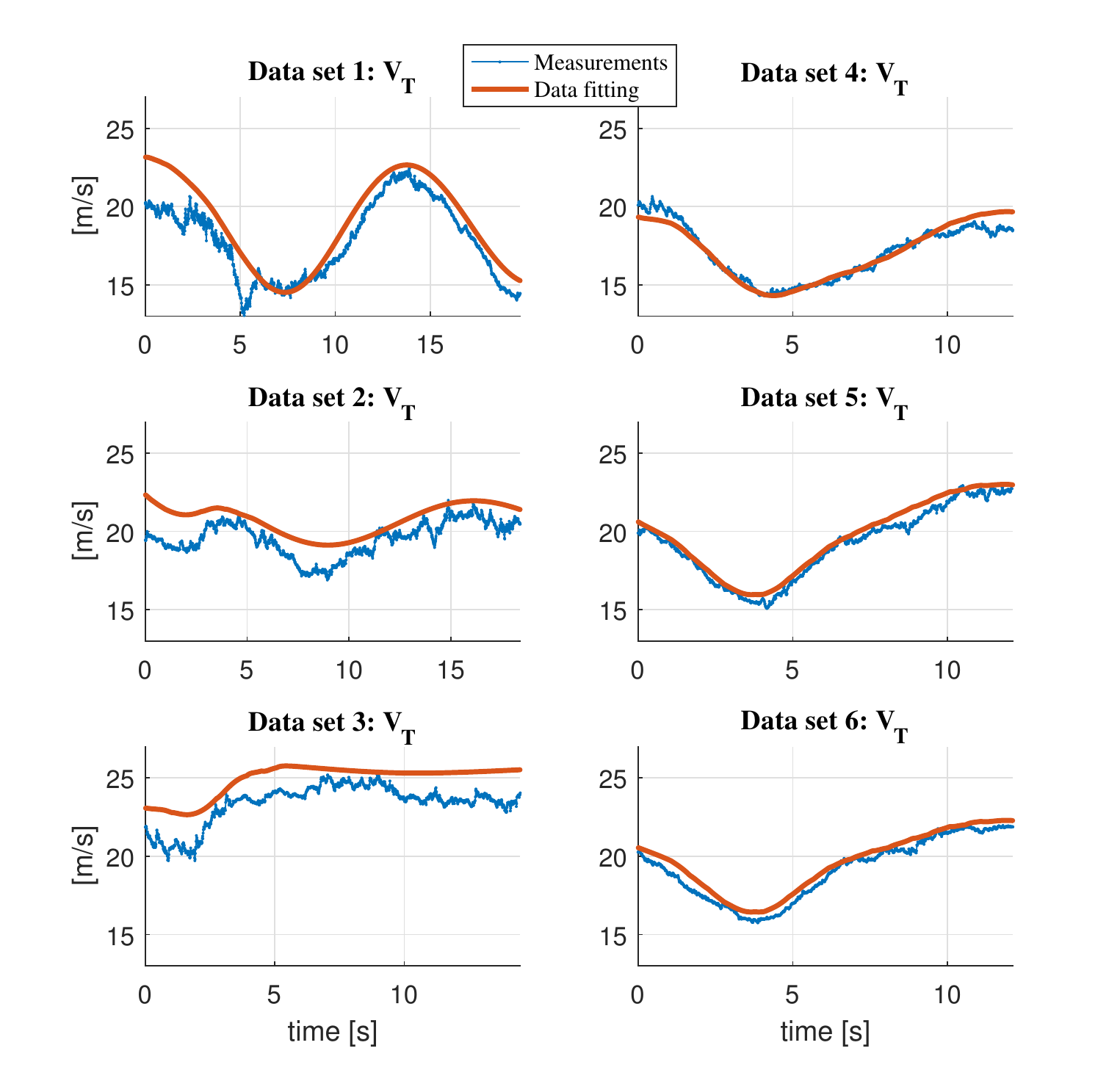} 
	\caption{Data fitting on multiple experiments along the longitudinal dynamics: airspeed $\VT$.}
	\label{fig:DataFitting_Airspeed}
\end{figure}

\begin{figure}[tbhp]
	\centering
	\includegraphics[scale=0.6]{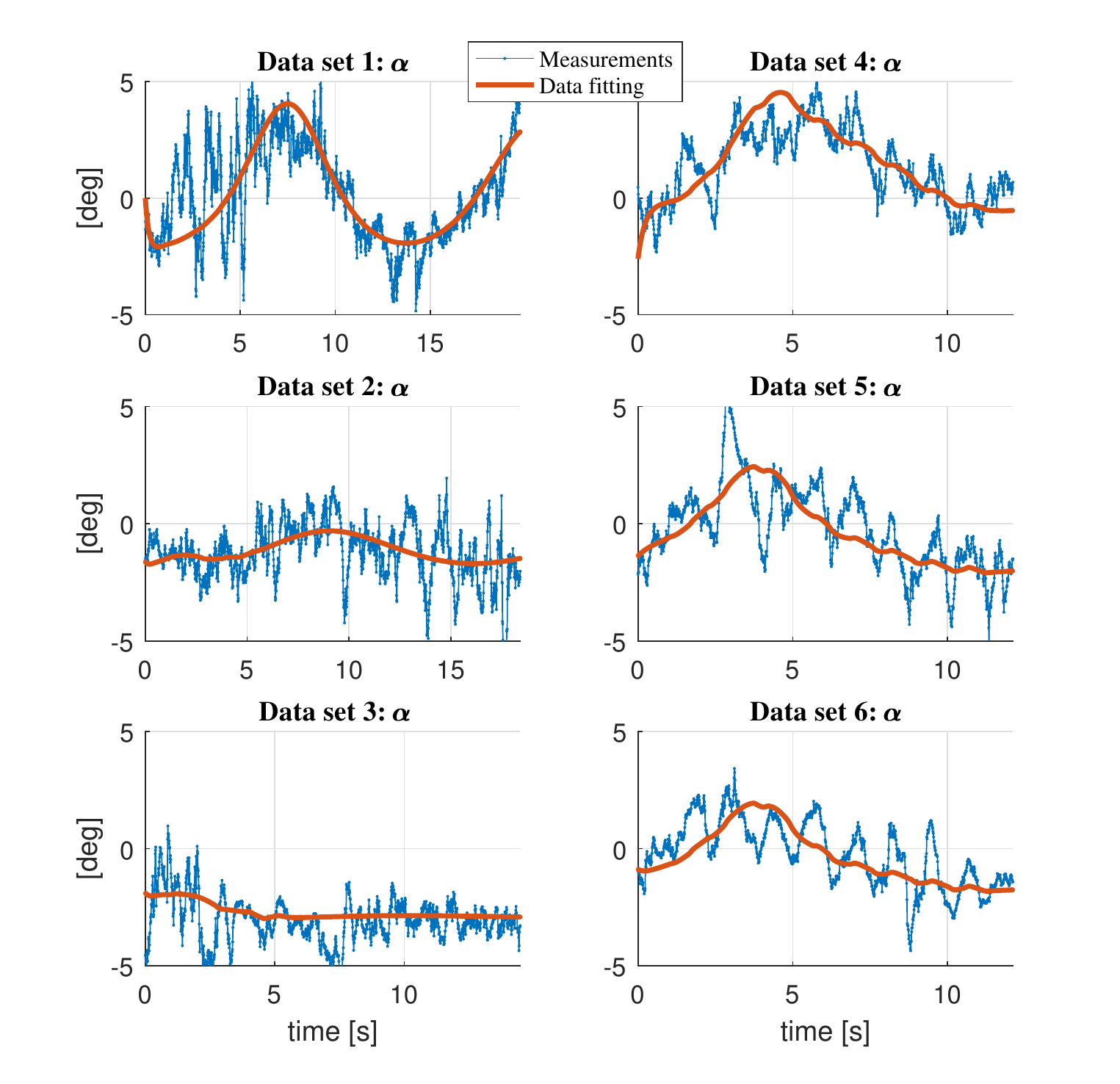} 
	\caption{Data fitting on multiple experiments along the longitudinal dynamics: angle of attack $\alpha$.}
	\label{fig:DataFitting_AngleOfAttack}
\end{figure}

\begin{figure}[tbhp]
	\centering
	\includegraphics[scale=0.6]{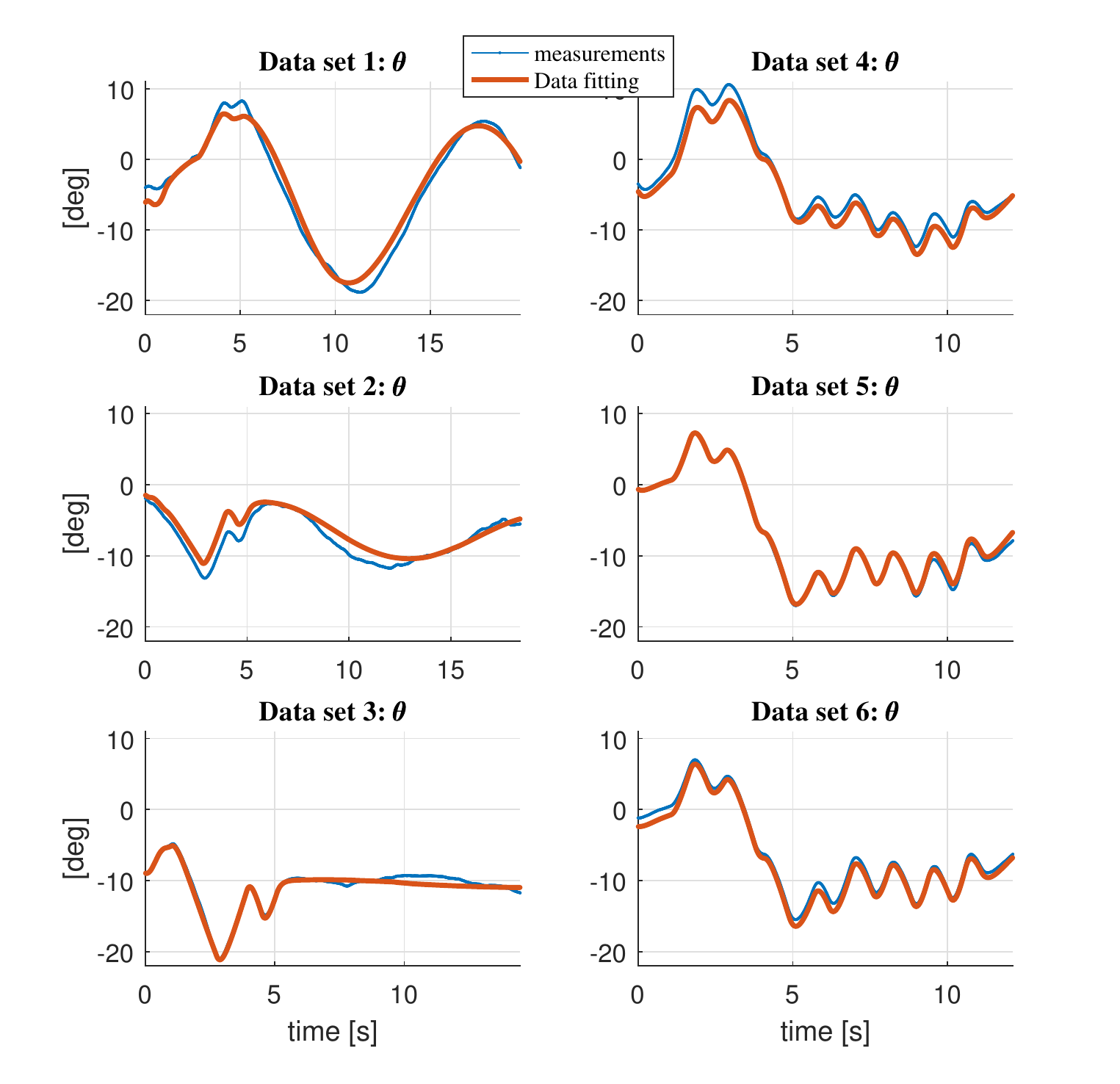} 
	\caption{Data fitting on multiple experiments along the longitudinal dynamics: pitch angle $\theta$.}
	\label{fig:DataFitting_PitchAngle}
\end{figure}

\begin{figure}[tbhp]
	\centering
	\includegraphics[scale=0.6]{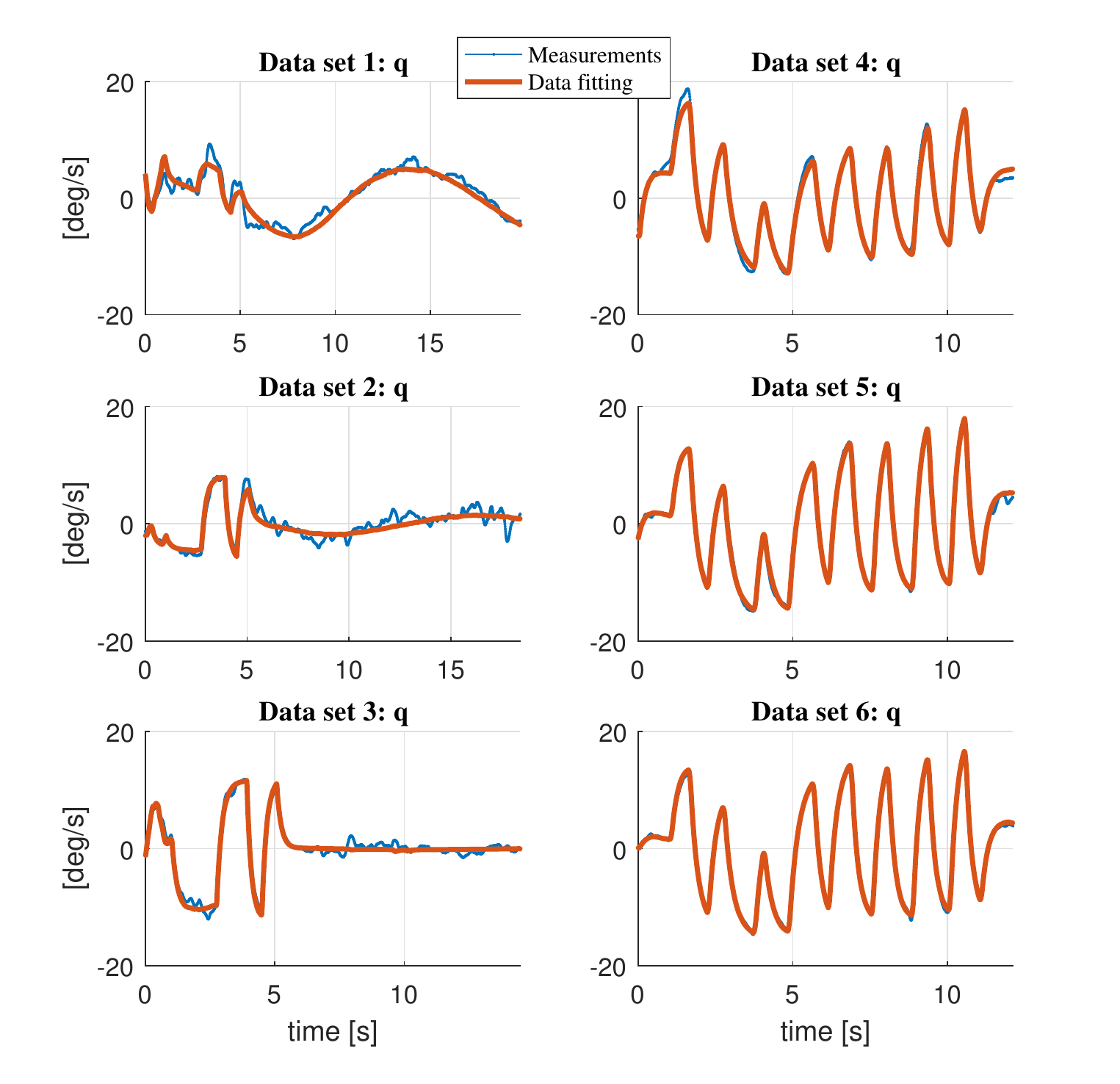} 
	\caption{Data fitting on multiple experiments along the longitudinal dynamics: pitch rate $q$.}
	\label{fig:DataFitting_PitchRate}
\end{figure}

Figures~\ref{fig:DataFitting_Airspeed},\ref{fig:DataFitting_AngleOfAttack},\ref{fig:DataFitting_PitchAngle} and \ref{fig:DataFitting_PitchRate} show the data fitting on the airspeed $\VT$, angle of attack $\alpha$, pitch angle $\theta$ and pitch rate $q$, respectively. Note that the measurements are suitably low-pass filtered using zero-lag filtering in order to focus on the rigid-body modes only. The control surface inputs are measured via feedback sensors on the aircraft, which allows the estimation to proceed without requiring knowledge of the actuator dynamics. The control surface deflection measurements have no discernible noise, though quantization errors equal to $0.25 \, ^\circ$ are present and compensated.\linebreak Finally, a one frame transport delay of the measurements is used. 

The overall data fitting is satisfactory except for the airspeed $\VT$, where biases arise in the conventional experiments where turbulence effect are more consistent. 

\subsection{Model validation}\label{sec:ModelValidation}
Because a significant inaccuracy on some derivatives relative to the Phugoid mode are expected (see section \ref{sec:PreliminaryAnalysis}) and biases on the airspeed data fitting are observed in section \ref{sec:dataFitting}, the estimates $C_{X_{q}}$, $C_{X_{\delta_{e}}}$ are set to their a priori values. 
As mentioned in section~\ref{sec:PreliminaryAnalysis}, in this way low frequency errors might arise in the identified model though, standard feedback controls can easily handle such model mismatch \cite{dorobantu2013system}.   

Furthermore, it turns out that the estimated derivative $C_{Z_{q}}$, i.e., the force variation along the $Z$-axis has no reasonable physical meaning and for this reason its value is fixed to the a priori estimate, too. At any way, uncertainties on $C_{Z_{q}}$ does not significantly deteriorate the predictive capability of the derived model \cite{mulder2000flight}.     

Table~\ref{tab:AerodynamicsDerivativesValidation} collects the a priori $\mathbf{\tilde{p}}$ and estimated $\mathbf{p}^{*}$ dimensionless aerodynamic longitudinal derivatives augmented with the set of parameters $\mathbf{p_{v}}$ which will be used for model validation. The set of parameters $\mathbf{p_{v}}$ provides the identified Phugoid and Short-period mode shown in Table~\ref{tab:EstimatedModes}.  

\begin{table}[tbhp]
	\caption{Dimensionless aerodynamic longitudinal derivatives}
	\label{tab:AerodynamicsDerivativesValidation}
	\centering
	\begin{tabular}{crrr}\hline
		$\mathbf{p}$        &  $\mathbf{\tilde{p}}$ & $\mathbf{p}^{*}$ & $\mathbf{p_{v}}$\\\hline
		$C_{X_{0}}$         &   -0.033      &    0.000     &    0.000   \\
		$C_{X_{\alpha}}$    &    0.409      &   -0.668     &    0.668   \\
		$C_{X_{q}}$			&   -0.603      &  -22.515     &   -0.603   \\
		$C_{X_{\delta_{e}}}$&   -0.011      &   -0.885     &   -0.011   \\
		$C_{Z_{0}}$			&   -0.528      &   -0.561     &   -0.561   \\
		$C_{Z_{\alpha}}$	&   -4.225      &   -5.012     &   -5.012   \\
		$C_{Z_{q}}$			&   -7.500      &  -61.940     &   -7.500   \\
		$C_{Z_{\delta_{e}}}$&   -0.310      &    0.122     &    0.122   \\
		$C_{m_{0}}$		    &   -0.031      &    0.061     &    0.061   \\
		$C_{m_{\alpha}}$ 	&   -0.607      &   -0.779     &   -0.779   \\
		$C_{m_{q}}$			&  -11.300      &  -24.923     &  -24.923   \\
		$C_{m_{\delta_{e}}}$&   -1.420      &   -1.004     &   -1.004   \\\hline		
	\end{tabular}
\end{table}
\begin{table}[tbhp]
	\caption{Identified longitudinal modes}
	\label{tab:EstimatedModes}
	\centering
	\begin{tabular}{cccc}
		\hline       
		\bf Mode       & \bf Short-period &\bf Phugoid &                 \\\hline
		$\omega_{n}$   &       5.548      &  0.587     & $\mathrm{rad/s}$\\
		$\tau $        &       0.180      &  1.704     & $\mathrm{s}$    \\
		$\delta$       &       0.843      &  0.036     & $-$    \\
		$S_{\%}$       &       0.721      & 89.210     & \%   \\
		$P_{O}$        &       2.108      & 10.712     & $\mathrm{s}$    \\\hline
	\end{tabular}
\end{table}

\begin{figure}[tbhp]
	\centering
	\includegraphics[scale=0.6]{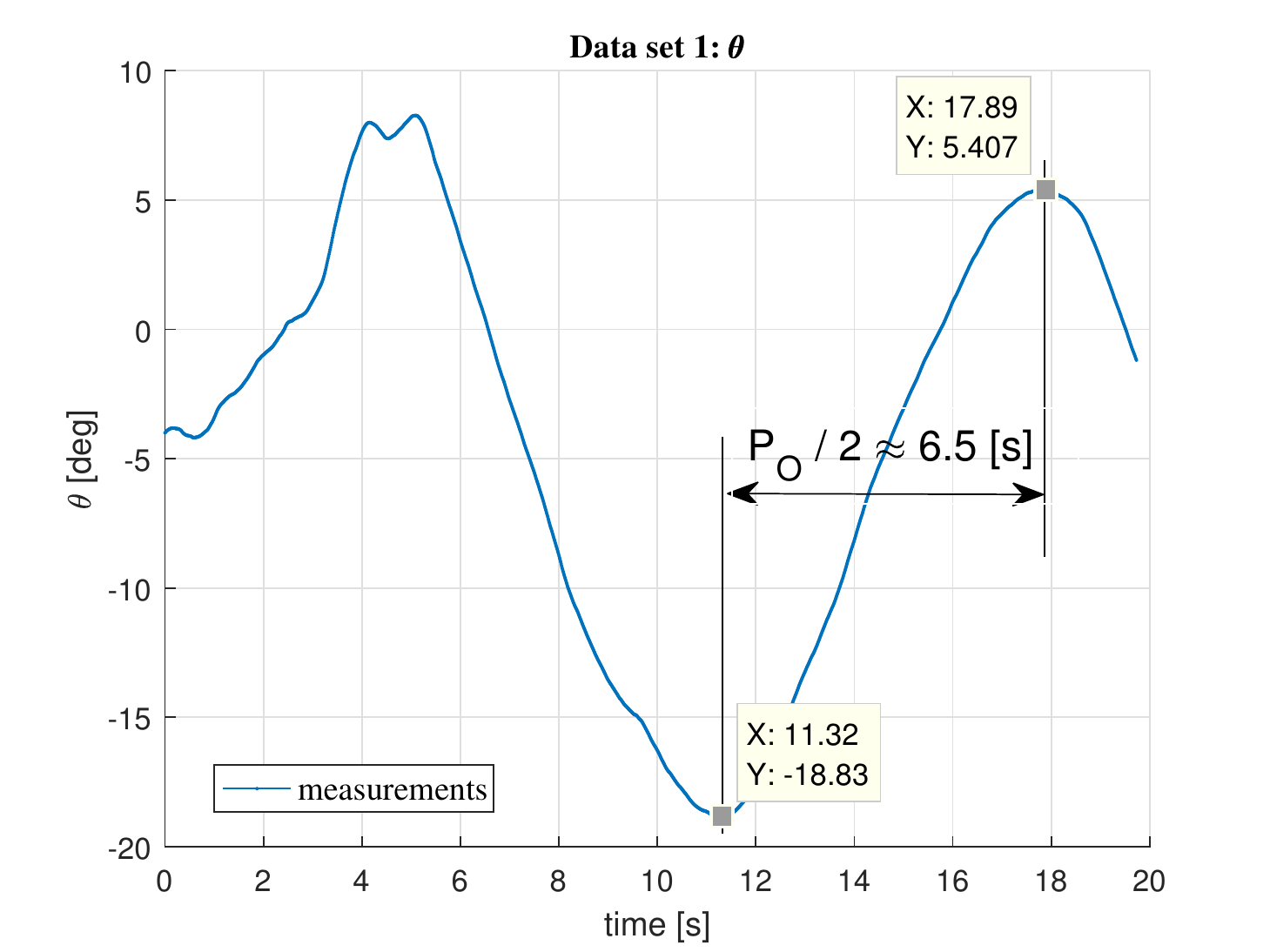} 
	\caption{Observed Phugoid period of oscillation $P_{\mathrm{O}}$.}
	\label{fig:EstimatedPhugoid_Po}
\end{figure}

Also in this case, a discrepancy is observed between the estimated Phugoid period ($P_{\mathrm{o}} \approx 11 \, \mathrm{s}$) and the observed one ($P_{\mathrm{o}} \approx 13 \, \mathrm{s}$) in the airspeed responses shown in figure~\ref{fig:EstimatedPhugoid_Po}.

\begin{figure}[tbhp]
	\centering
	\includegraphics[scale=0.6]{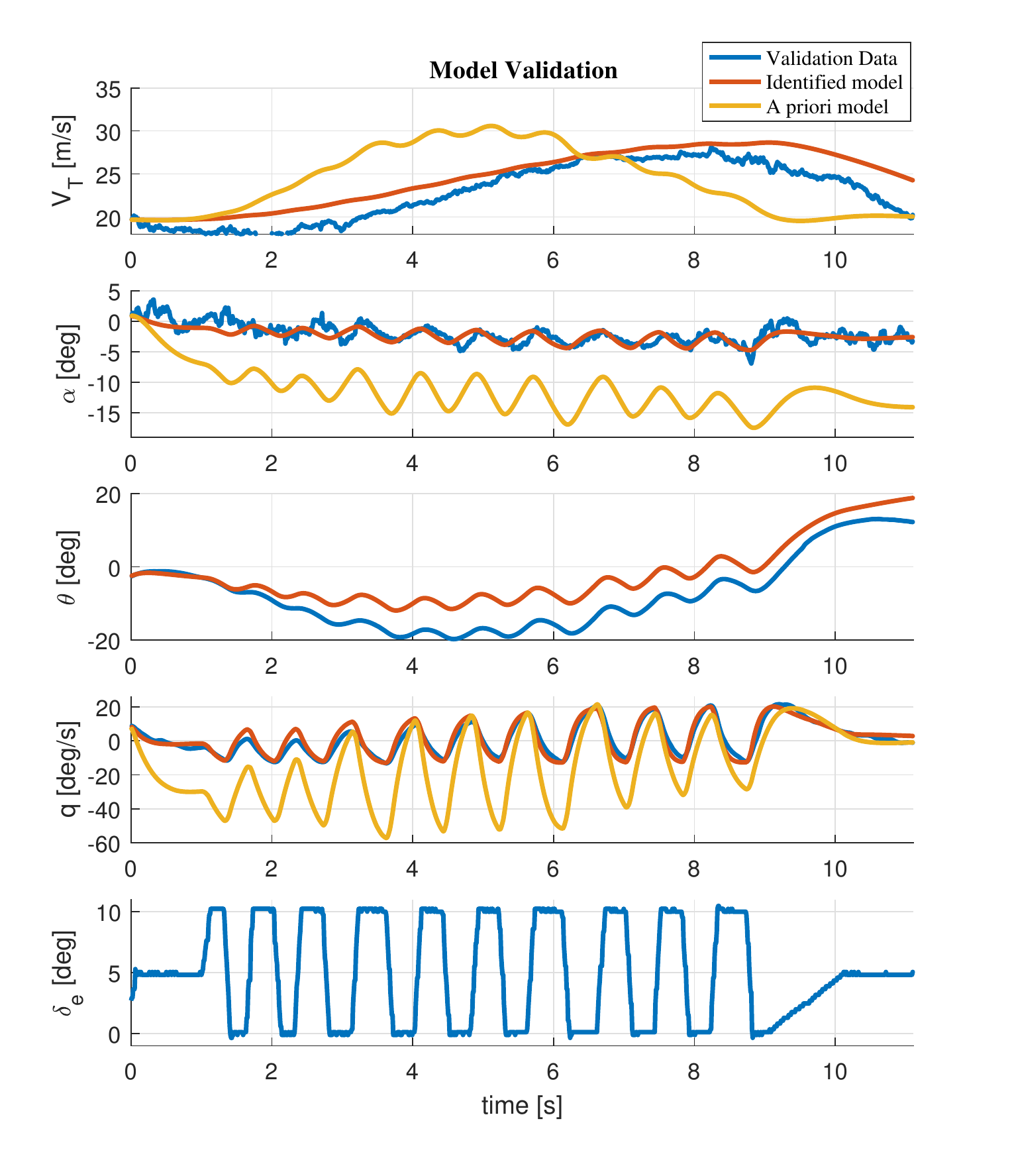} 
	\caption{Model structure assessment via validation data set. The a priori pitch angle $\theta$ response is not shown due to its large deviation w.r.t. the obtained experimental values.}
	\label{fig:ModelValidation}
\end{figure}

\begin{figure}[tbhp]
	\centering
	\includegraphics[scale=0.6]{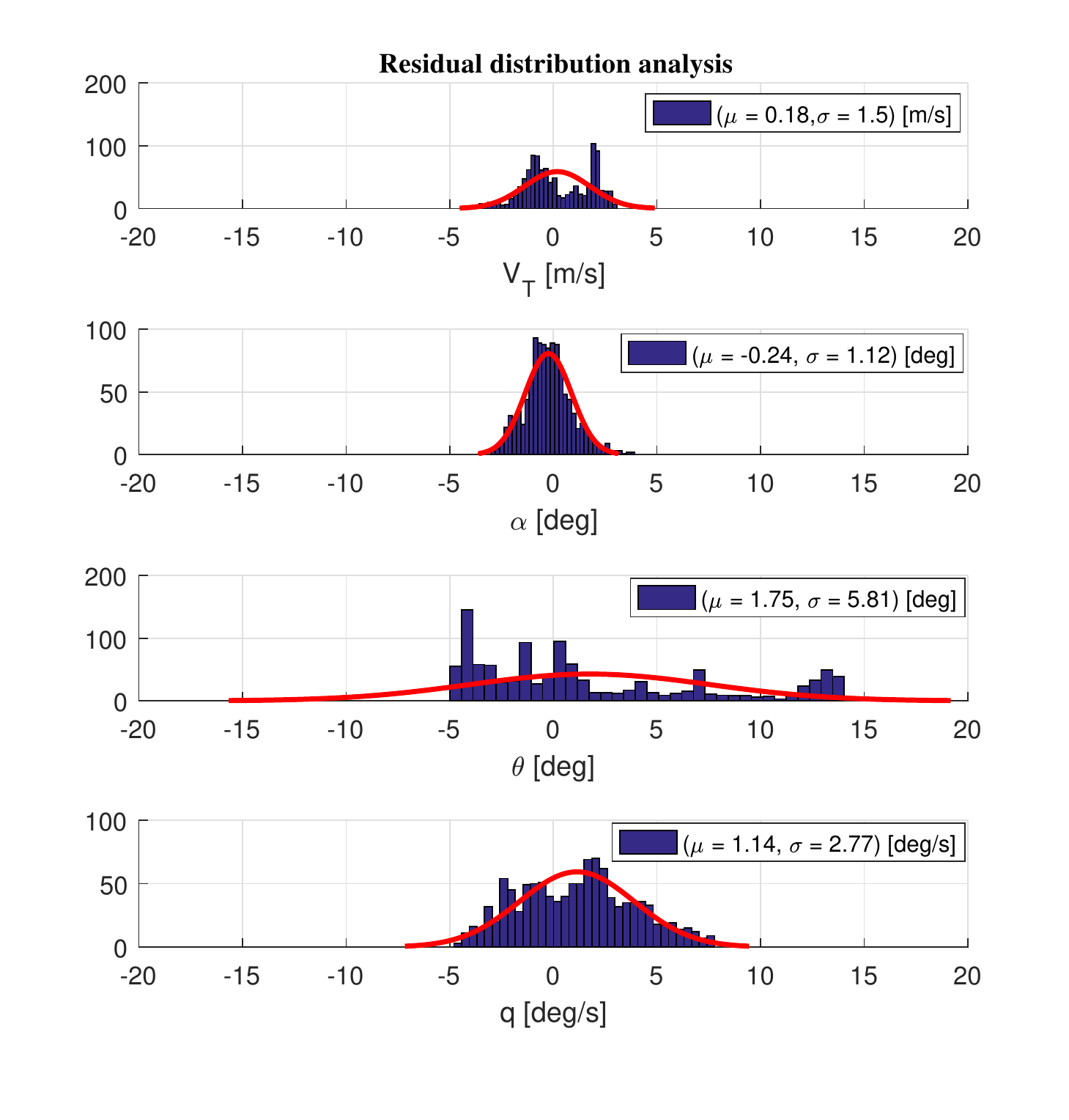} 
	\caption{Residual distribution analysis on validation data set with corresponding mean value $\mu$ and standard deviation $\sigma$.}
	\label{fig:ResidualDistribution}
\end{figure}

The accuracy of an identified model is ultimately assessed via its capability to predict time responses \cite{dorobantu2013system}. For validation purpose, the identified model is simulated using a further flight test experiment shown in figure~\ref{fig:ModelValidation}. One can observe that the identified model provides a better fitting compared to the a priori one despite inaccuracies on the Phugoid mode. Figure~\ref{fig:ResidualDistribution} shows the corresponding residual distributions $\epsilon$ defined as
\begin{equation}\label{eq:residual_analysis}
\mathbf{\epsilon_{k}} = \:  \hat{\mathbf{y}}_{k} - \mathbf{h}\left (\mathbf{x}_{k},\hat{\mathbf{u}}_{k},\mathbf{p_{v}} \right), k = 1,\ldots,N_{\mathrm{v}} \\
\end{equation}
with $N_{\mathrm{v}}$ the number of samples related to the validation data set. 
Practically speaking, the residual is the part of the data that the model is not able to reproduce; the aim is to achieve a residual resembling a white noise signal. However, it is well-known that the residuals will not be white noise if the real system has significant process noise (atmospheric turbulence) \cite{licitra2017pe}. 

%In this case the airspeed and pitch angle residual differ consistently from a Gaussian distribution and since these dynamics one can observer that process noise provide in general bias on the estimate and in this case 

Finally, estimation results are assessed via the so called \ac{TIC} which is defined by the following relationship \cite{cook2007flight}
\begin{equation}\label{eq:TIC}
\mathrm{TIC} = \frac{\sqrt{\frac{1}{N_{\mathrm{v}}}  \Sigma_{i=1}^{N_{\mathrm{v}}} \left( \mathbf{\hat{y}}_{i} - \mathbf{h}\left (\mathbf{x}_{i},\hat{\mathbf{u}}_{i},\mathbf{p_{v}} \right) \right)^{2} }}{\sqrt{\frac{1}{N_{\mathrm{v}}}  \Sigma_{i=1}^{N_{\mathrm{v}}} \mathbf{\hat{y}}_{i}^{2}} + \sqrt{\frac{1}{N_{\mathrm{v}}}  \Sigma_{i=1}^{N_{\mathrm{v}}}  \mathbf{h}\left (\mathbf{x}_{i},\hat{\mathbf{u}}_{i},\mathbf{p_{v}} \right)^{2} }}
\end{equation}

The \ac{TIC} provides a basis of judgment regarding the degree of predictability of an mathematical (estimated) model via a normalized metric between 0 and 1. A value of $\mathrm{TIC} = 0$ denotes a perfect match whereas and $\mathrm{TIC} = 1$ indicates the worst case scenario i.e. the mathematical model is not able to explain any of the data. Values of $\ac{TIC} \le 0.25$ correspond to accurate prediction for rigid wing aircraft \cite{tischler2006aircraft,jategaonkar2004aerodynamic}. Table~\ref{tab:TIC} summarizes the \ac{TIC} values for this work. 

\begin{table}[tbhp]
	\caption{Theil Inequality Coefficients}
	\label{tab:TIC}
	\centering
	\begin{tabular}{ccccc}
		\hline       
		& $\VT$ & $\alpha$ & $\theta$ & $q$ \\\hline
		TIC   & 0.04  &  0.22    &  0.29    & 0.15 \\\hline
	\end{tabular}
\end{table}

Results shows that the pitch rate $q$ is captured with high accuracy as well as the airspeed response $V_{\mathrm{T}}$ despite the uncertainties mentioned above whereas the angle of attack $\alpha$ is near to the limit for accurate models. Finally, the pitch angle $\theta$ is slightly over the limit. These model mismatch might be attributed to 
\begin{itemize}
	\item Difficulties in obtaining good excitation on the angle of attack response without flight envelope violation and with sufficient \ac{SNR};
	\item low uncertainties on the damping dynamics in the Short-period mode can causes unsatisfactory pitch angle prediction due to integration errors;
	\item neglected coupling between the longitudinal and lateral motion.
\end{itemize}
One has to point out that a typical control strategy of a rigid wing \ac{AWES} is based on the angle of attack tracking during \textit{power production} phase \cite{Cap26AWEbook}. This results suggest to include on such control strategy an integrator term in the forward path relative to the angle of attack channel in order to account for possible model mismatches, disturbances and at the same time achieve reasonable tracking performance.  

\section{Conclusions}\label{sec:Conclusions}
In this paper, real flight test experiments and subsequent time domain Model-Based Parameter Estimations (MBPEs) have been carried out for a rigid wing Airborne Wind Energy System (AWES). A suitable and comprehensive non-linear mathematical model for (MBPE) problems relative to the aircraft dynamics was introduced and an overview of the flight test procedure has been provided. The experimental data were obtained for the longitudinal dynamics for the steady state wing-level trim condition. The optimization problem was initialized using a priori aerodynamic derivatives obtained via the lifting line method. Finally, the identified model was assessed by time domain model validation, residual distribution analysis and Theil Inequality Coefficients (TIC). 

Experimental results have shown that system identification via real flight tests is able to improve the predictive capability of low fidelity a priori models for rigid wing (AWES). However, baseline models are equally important to deal with non identifiable dynamics as well as for designing maneuvers for system identification purposes. 

Future work will aim to the implementation of parameter estimation algorithms which are robust with respect to turbulence e.g. Filter-Error method \cite{morelli2015filter}.

\section*{Acknowledgments}
The authors would like to thank Ampyx Power B.V., and in particular the flight operations team for conducting the system identification flight tests.

%\section*{References}
\bibliography{mybibfile}

\begin{thebibliography}{10}
\expandafter\ifx\csname url\endcsname\relax
  \def\url#1{\texttt{#1}}\fi
\expandafter\ifx\csname urlprefix\endcsname\relax\def\urlprefix{URL }\fi
\expandafter\ifx\csname href\endcsname\relax
  \def\href#1#2{#2} \def\path#1{#1}\fi

\bibitem{loyd1980crosswind}
M.~L. Loyd, Crosswind kite power (for large-scale wind power production),
  Journal of energy 4~(3) (1980) 106--111.

\bibitem{AP}
Ampyx power: Airborne wind energy, \url{https://www.ampyxpower.com} (2017).

\bibitem{makani}
Makani, \url{https://x.company/makani/} (2017).

\bibitem{twingtec}
Twingtec: Wind energy 2.0, \url{http://twingtec.ch/} (2017).

\bibitem{kitemill}
kitemill, \url{http://www.kitemill.com/} (2017).

\bibitem{e-kite}
e-kite, \url{http://www.e-kite.com/} (2017).

\bibitem{anderson2017fundamentals}
J.~D. Anderson~Jr, Fundamentals of aerodynamics, 6th Edition, McGraw-Hill,
  2017.

\bibitem{hoak1975usaf}
D.~Hoak, R.~Finck, USAF stability and control DATCOM, National Technical
  Information Service, 1975.

\bibitem{licitra2017pe}
G.~Licitra, P.~Williams, J.~Gillis, S.~Ghandchi, S.~Sieberling, R.~Ruiterkamp,
  M.~Diehl, Aerodynamic parameter identification for an airborne wind energy
  pumping system, in: IFAC 2017 World Congress, Toulouse, France. 9-14 July,
  2017.

\bibitem{diehl2013airborne}
M.~Diehl, Airborne wind energy: Basic concepts and physical foundations, in:
  Airborne Wind Energy, Springer, 2013, pp. 3--22.
\newblock \href {http://dx.doi.org/10.1007/978-3-642-39965-7_1}
  {\path{doi:10.1007/978-3-642-39965-7_1}}.

\bibitem{Cap26AWEbook}
R.~Ruiterkamp, S.~Sieberling, Description and preliminary test results of a six
  degrees of freedom rigid wing pumping system, in: Airborne wind energy,
  Springer, 2013, pp. 443--458.
\newblock \href {http://dx.doi.org/10.1007/978-3-642-39965-7_26}
  {\path{doi:10.1007/978-3-642-39965-7_26}}.

\bibitem{PowerPatternvideo}
{Ampyx Power B.V., Working Principle},
  \url{https://www.youtube.com/watch?v=32CYqjoRnYY}, [Online; accessed
  06-November-2017] (2017).

\bibitem{licitra2017oed}
G.~Licitra, A.~B\"{u}rger, P.~Williams, R.~Ruiterkamp, M.~Diehl, Optimum
  experimental design of a rigid wing awe pumping system, in: IEEE 56th
  Conference on Decision and Control, 2017.

\bibitem{licitra2017optinput}
G.~Licitra, A.~B\"{u}rger, P.~Williams, R.~Ruiterkamp, M.~Diehl, Optimal input
  design for autonomous aircraft, submitted to Journal of Control Engineering
  Practice.

\bibitem{gros2013modeling}
S.~Gros, M.~Diehl, Modeling of airborne wind energy systems in natural
  coordinates, in: Airborne wind energy, Springer, 2013, pp. 181--203.
\newblock \href {http://dx.doi.org/10.1007/978-3-642-30023-327}
  {\path{doi:10.1007/978-3-642-30023-327}}.

\bibitem{licitra2016optimal}
G.~Licitra, S.~Sieberling, S.~Engelen, P.~Williams, R.~Ruiterkamp, M.~Diehl,
  Optimal control for minimizing power consumption during holding patterns for
  airborne wind energy pumping system, in: Control Conference (ECC), 2016
  European, IEEE, 2016, pp. 1574--1579.
\newblock \href {http://dx.doi.org/10.1109/ECC.2016.7810515}
  {\path{doi:10.1109/ECC.2016.7810515}}.

\bibitem{licitra2017viability}
G.~Licitra, J.~Koenemann, G.~Horn, P.~Williams, R.~Ruiterkamp, M.~Diehl,
  Viability assessment of a rigid wing airborne wind energy pumping system, in:
  Process Control (PC), 2017 21st International Conference on, IEEE, 2017, pp.
  452--458.
\newblock \href {http://dx.doi.org/10.1109/PC.2017.7976256}
  {\path{doi:10.1109/PC.2017.7976256}}.

\bibitem{williams2007modeling}
P.~Williams, B.~Lansdorp, W.~Ockels, in: Modeling and control of a kite on a
  variable length flexible inelastic tether, AIAA Guidance, navigation and
  control conference, 2007.

\bibitem{williams2008modeling}
P.~Williams, B.~Lansdorp, R.~Ruiterkamp, W.~Ockels, in: Modeling, simulation,
  and testing of surf kites for power generation, AIAA Modeling and Simulation
  Technologies Conference and Exhibit, 2008.
\newblock \href {http://dx.doi.org/10.2514/6.2008-6693}
  {\path{doi:10.2514/6.2008-6693}}.

\bibitem{AWEbook}
M.~Diehl, R.~Schmehl, U.~Ahrens, Airborne Wind Energy. Green Energy and
  Technology, Springer Berlin Heidelberg, 2014.
\newblock \href {http://dx.doi.org/10.1007/978-3-642-39965-7}
  {\path{doi:10.1007/978-3-642-39965-7}}.

\bibitem{stevens2015aircraft}
B.~L. Stevens, F.~L. Lewis, E.~N. Johnson, Aircraft control and simulation:
  dynamics, controls design, and autonomous systems, Wiley, 2015.
\newblock \href {http://dx.doi.org/10.1002/9781119174882}
  {\path{doi:10.1002/9781119174882}}.

\bibitem{mulder2000flight}
J.~Mulder, W.~Van~Staveren, J.~van~der Vaart, Flight dynamics (lecture notes),
  TU Delft, 2000.

\bibitem{etkin2012dynamics}
B.~Etkin, Dynamics of atmospheric flight, Wiley, 1972.

\bibitem{goman1994state}
M.~Goman, A.~Khrabrov, State-space representation of aerodynamic
  characteristics of an aircraft at high angles of attack, Journal of Aircraft
  31~(5) (1994) 1109--1115.
\newblock \href {http://dx.doi.org/10.2514/3.46618}
  {\path{doi:10.2514/3.46618}}.

\bibitem{de1987accurate}
R.~De~Jong, J.~Mulder, Accurate estimation of aircraft inertia characteristics
  from a single suspension experiment, Journal of Aircraft 24~(6) (1987)
  362--370.
\newblock \href {http://dx.doi.org/10.2514/3.45454}
  {\path{doi:10.2514/3.45454}}.

\bibitem{lyons2002obtaining}
D.~Lyons, Obtaining optimal results with filar pendulums for moment of inertia
  measurements, Society of Allied Weight Engineers 62~(2) (2002) 5--20.

\bibitem{soderstrom2001identification}
T.~Soderstrom, P.~Stoica, System Identification, Prentice-Hall, 2001.
\newblock \href {http://dx.doi.org/10.1007/BF01211647}
  {\path{doi:10.1007/BF01211647}}.

\bibitem{dorobantu2013system}
A.~Dorobantu, A.~Murch, B.~Mettler, G.~Balas, System identification for small,
  low-cost, fixed-wing unmanned aircraft, Journal of Aircraft 50~(4) (2013)
  1117--1130.
\newblock \href {http://dx.doi.org/10.2514/1.C032065}
  {\path{doi:10.2514/1.C032065}}.

\bibitem{klein1990optimal}
E.~Morelli, V.~Klein, Optimal input design for aircraft parameter estimation
  using dynamic programming principles, in: AIAA Atmospheric Flight Mechanics
  Conference, 1990.

\bibitem{cook2007flight}
M.~Cook, Flight Dynamics Principles: A Linear Systems Approach to Aircraft
  Stability and Control, 3rd Edition, Elsevier, 2013.
\newblock \href {http://dx.doi.org/10.1017/S0001924000010873}
  {\path{doi:10.1017/S0001924000010873}}.

\bibitem{mcruer2014aircraft}
D.~T. McRuer, D.~Graham, I.~Ashkenas, Aircraft dynamics and automatic control,
  Princeton University Press, 2014.

\bibitem{koehler1977auslegung}
R.~Koehler, K.~Wilhelm, Auslegung von eingangssignalen f{\"u}r die
  kennwertermittlung, DFVLR-IB (1977) 154--77.

\bibitem{mulder1994identification}
J.~Mulder, J.~Sridhar, J.~Breeman, Identification of dynamic system -
  application to aircraft. Part 2: nonlinear analysis and manoeuvre design,
  Advisory Group for Aerospace Research \& Development (AGARD), 1994.

\bibitem{marchand1977untersuchung}
M.~Marchand, Untersuchung der bestimmbarkeit der flugmechanischen derivative
  des ccv-versuchstr{\"a}gers f-104 g, Tech. rep. (1977).

\bibitem{OptSysIDvideo}
{Ampyx Power B.V., Optimal System Identification Flight Test},
  \url{https://www.youtube.com/watch?v=KBq5TTTOqf8}, [Online; accessed
  17-September-2017] (2017).

\bibitem{versteeg2007introduction}
H.~K. Versteeg, W.~Malalasekera, An introduction to computational fluid
  dynamics - the finite volume method, Addison-Wesley-Longman, 1995.

\bibitem{morelli1993oed}
E.~A. Morelli, Practical input optimization for aircraft parameter estimation
  experiments, National Aeronautics and Space Administration, Langley Research
  Center, 1993.

\bibitem{tischler2006aircraft}
M.~B. Tischler, R.~K. Remple, Aircraft and rotorcraft system identification,
  AIAA education series, 2006.
\newblock \href {http://dx.doi.org/10.2514/4.868207}
  {\path{doi:10.2514/4.868207}}.

\bibitem{morelli2006practical}
E.~A. Morelli, Practical aspects of the equation-error method for aircraft
  parameter estimation, in: AIAA Atmospheric Flight Mechanics Conference, Vol.
  6144, 2006.

\bibitem{ljung1998system}
L.~Ljung, System identification: Theory for the User, Prentice Hall PTR, 1998.
\newblock \href {http://dx.doi.org/10.1002/047134608X.W1046}
  {\path{doi:10.1002/047134608X.W1046}}.

\bibitem{diehl2014ocp}
M.~Diehl, Optimal Control and Estimation (lecture notes), University of
  Freiburg, 2014.

\bibitem{bock1984multiple}
H.~G. Bock, K.-J. Plitt, in: A multiple shooting algorithm for direct solution
  of optimal control problems, Proceedings of the IFAC world congress, 1984.
\newblock \href {http://dx.doi.org/10.1016/S1474-6670(17)61205-9}
  {\path{doi:10.1016/S1474-6670(17)61205-9}}.

\bibitem{bock2013model}
H.~G. Bock, T.~Carraro, W.~J{\"a}ger, S.~K{\"o}rkel, R.~Rannacher,
  J.~Schl{\"o}der, Model Based Parameter Estimation: Theory and Applications,
  Vol.~4, Springer Science \& Business Media, 2013.
\newblock \href {http://dx.doi.org/10.1007/978-3-642-30367-8}
  {\path{doi:10.1007/978-3-642-30367-8}}.

\bibitem{Andersson2013b}
J.~Andersson, {A} {G}eneral-{P}urpose {S}oftware {F}ramework for {D}ynamic
  {O}ptimization, Ph.D. thesis, Department of Electrical Engineering (ESAT/SCD)
  and Optimization in Engineering Center, KU Leuven (2013).

\bibitem{Waechter2006}
A.~W{\"a}chter, L.~T. Biegler, On the implementation of an interior-point
  filter line-search algorithm for large-scale nonlinear programming,
  Mathematical Programming 106~(1) (2006) 25--57.
\newblock \href {http://dx.doi.org/10.1007/s10107-004-0559-y}
  {\path{doi:10.1007/s10107-004-0559-y}}.

\bibitem{HSL2017}
HSL, \href{http://www.hsl.rl.ac.uk/}{A collection of fortran codes for large
  scale scientific computation} (2017).
\newline\urlprefix\url{http://www.hsl.rl.ac.uk/}

\bibitem{jategaonkar2004aerodynamic}
R.~Jategaonkar, D.~Fischenberg, W.~von Gruenhagen, Aerodynamic modeling and
  system identification from flight data-recent applications at dlr, Journal of
  Aircraft 41~(4) (2004) 681--691.
\newblock \href {http://dx.doi.org/10.2514/1.3165} {\path{doi:10.2514/1.3165}}.

\bibitem{morelli2015filter}
G.~J. A., E.~A. Morelli, A new formulation of the filter-error method for
  aerodynamic parameter estimation in turbulence, in: AIAA Atmospheric Flight
  Mechanics Conference, 2015.
\newblock \href {http://dx.doi.org/10.2514/6.2015-2704}
  {\path{doi:10.2514/6.2015-2704}}.

\end{thebibliography}

\end{document}